\documentclass[10.5pt,reqno]{amsart}
\usepackage{amsmath,amsfonts,amssymb,amsthm,mathtools}
\usepackage[numbers,sort&compress]{natbib}
\usepackage{mathrsfs}
\usepackage{microtype}
\usepackage{enumitem}
\usepackage[colorlinks=true]{hyperref}
\newtheorem{theorem}{Theorem}[section]
\newtheorem{lemma}[theorem]{Lemma}
\newtheorem{proposition}[theorem]{Proposition}
\newtheorem{corollary}[theorem]{Corollary}

\theoremstyle{definition}
\newtheorem{definition}[theorem]{Definition}
\theoremstyle{remark}
\newtheorem{remark}[theorem]{Remark}
\numberwithin{equation}{section}

\newcommand{\E}{\mathbb{E}}
\newcommand{\T}{\mathbb{T}}
\newcommand{\R}{\mathbb{R}}
\newcommand{\C}{\mathbb{C}}
\newcommand{\dd}{\,\mathrm{d}}
\newcommand{\norm}[1]{\left\|#1\right\|}
\newcommand{\abs}[1]{\left|#1\right|}

\begin{document}
\title{Null and approximate controllability of stochastic parabolic equations from measurable space--time sets}
\author{\sffamily Yuanhang Liu$^{1^*}$ and Weijia Wu$^{2}$ \\
{\sffamily\small $^1$ School of Mathematics and Statistics, Zhaoqing University, Guangdong 526061, China.}\\
{\sffamily\small $^2$ School of Mathematics and Physics, North China Electric Power University, Beijing, 102206, China. }\\
}
\thanks{$ ^*$ Corresponding author: liuyuanhang97@163.com}
\email{liuyuanhang97@163.com}
\email{weijiawu@yeah.net}

\begin{abstract}
In this paper, we consider a stochastic parabolic equation on the torus with drift, a space-dependent analytic noise coefficient and random initial data. A single control acts on the drift and is supported in an arbitrary measurable subset of the space-time cylinder, not necessarily a product set. The analytic norms of the coefficients need only be square integrable in time. A complex translation of the Fourier-Galerkin system yields quantitative analyticity of the adjoint state. Combining it with a Hilbert-valued propagation of smallness and a telescoping argument over level sets of the observation measure gives an observability estimate whose constant depends on the observation set only through a lower bound on its measure. We deduce null controllability with one adapted control and characterize approximate controllability by the accumulation of the observation set at the terminal time. Null controllability therefore does not imply approximate controllability, the obstruction being the filtration structure of the target.
\end{abstract}
\maketitle
\raggedbottom
\pagestyle{myheadings}
\thispagestyle{plain}
\markboth{Controllability from measurable space--time sets}{Y. Liu and W. Wu}


\section{Introduction}\label{sec:intro}

\subsection*{Motivation and deterministic background}

The controllability of parabolic equations is a classical chapter of control theory. In the deterministic setting the duality method reduces null controllability to an observability estimate for the adjoint equation, and the latter can be obtained from global Carleman estimates or from spectral inequalities combined with a time iteration. The situation changes when the equation is subject to multiplicative noise. The adjoint system is then a backward stochastic parabolic equation, which carries both a state and a martingale integrand. The integrand enters the drift of the adjoint equation and behaves as a nonhomogeneous term. This is the source of the main difficulty. It explains why the Carleman-based controllability results for forward stochastic parabolic equations use two controls, one acting on the drift and one acting on the diffusion \cite{TZ,HLP}, while single-control results were available only under additional structure on the noise coefficient, for instance when it does not depend on the space variable \cite{Lu2011,YZ,LWYZ2026}. This paper addresses two questions. The first is whether one control, acting only on the drift, can suffice for a genuinely space-dependent noise coefficient. The second concerns the geometry of the control region, which may be merely measurable with positive measure and may be a genuine subset of the space--time cylinder rather than a product of a spatial set and a time interval. Both questions are answered below for a class of stochastic parabolic equations on the torus, with space-dependent analytic coefficients, random initial data and a single drift control.

For the heat equation, observability from sets of positive measure is well understood. The classical spectral inequality \cite{LebeauRobbiano} gives observability from open spatial sets over the whole time interval, and global Carleman estimates give the same conclusion for parabolic equations with variable coefficients \cite{FursikovImanuvilov}. Measurable spatial sets over the whole time interval were treated in \cite{AE2013} for self-adjoint elliptic operators with analytic coefficients near the observation set, and measurable space--time sets contained in a fixed cylinder were treated in \cite{AEW}. When the observation set is measurable only in time, a partition adapted to its density points replaces the geometric partition of the time axis \cite{PhungWang2013}. Bounded space--time potentials and arbitrary measurable space--time sets were treated in \cite{WZ26}, where the telescoping is organized by the level sets of the accumulated measure of the observation times and the resulting subquadratic criterion is isolated as a separate combinatorial statement. One-dimensional problems with bounded space--time potentials are treated in \cite{Zhu2026}, and optimal constants for time-independent bounded measurable potentials in one dimension are obtained in \cite{ZhuZhuge2025}. The results based on measurable observation sets share a common structure, namely a one-time interpolation inequality with an exponential cost in the inverse time gap, an energy propagation and a telescoping argument over a partition of the time axis. The deterministic theory thus supplies the model for the interpolation and telescoping steps below.

\subsection*{Stochastic observability and the analytic coefficient class}

For stochastic parabolic equations the adjoint system is a backward stochastic equation carrying both a state and a martingale integrand. The integrand behaves as a nonhomogeneous term and destroys the usual Carleman scheme. Carleman-based controllability results are available for fourth-order and semilinear equations \cite{LuWang2022,ZhangXuLiu2026}. Quantitative unique continuation for the adjoint state, which is the analytic core of the measurable-set approach, requires separate arguments \cite{LY,LYZ2024}. These results yield unique continuation and, in some cases, observability from an open spatial region and a measurable time set, but they do not give observability from a spatial set of positive measure with empty interior.

A different route is available when the noise coefficient is independent of the space variable. Multiplication by such a coefficient commutes with the spectral projectors of the spatial operator, the projected adjoint equations close, and a spectral inequality combined with time iteration gives observability, hence null controllability, with a single drift control. This was first carried out for observation sets that are open in space and measurable in time \cite{Lu2011}, and later for measurable spatial sets as well \cite{YZ}. An abstract version covering several evolution equations is given in \cite{LWYZ2026}. In these results the noise coefficient is scalar in space, and this restriction is essential for the spectral argument. Approximate controllability requires in addition that the observation times accumulate at the terminal time, a condition identified in \cite{Lu2011} for a scalar noise coefficient. A spatially varying coefficient mixes frequencies. The unobserved high-frequency part of the martingale integrand can enter the low-frequency equation for the state, so the projected equations no longer close. This is the obstruction that the present paper has to overcome. A computation on $\T^{1}$ with a single-mode martingale integrand illustrating this mechanism is given in Remark \ref{rem:mixing}.

A recent work \cite{LW26} identified a class of spatially varying noise coefficients for which the frequency mixing can be controlled. The spatial profiles admit a bounded holomorphic extension to a complex polystrip of fixed width, and the analytic norm is only required to be square integrable in time, uniformly along sample paths. A direct analytic energy estimate based on a time-dependent complex translation of the Fourier--Galerkin equations yields quantitative spatial analyticity of the adjoint state without spatial regularity of the martingale integrand. Propagation of smallness in the Hilbert-valued setting then gives a one-time interpolation from any measurable spatial set of positive measure, and a telescoping argument gives observability from a cylindrical observation set, with an $L^{1}$-in-time mean square norm. Null controllability with a single adapted drift control and approximate controllability to arbitrary square integrable random targets follow. In \cite{LW26} the observation set is a cylinder, the drift coefficient is absent and the initial datum is deterministic. Approximate controllability is proved but not characterized.

\subsection*{Contributions and comparison with existing work}

The results of \cite{LW26} leave several questions open. The observation set is a cylinder, the equation has no drift, the initial datum is deterministic, and the approximate controllability statement is sufficient but not necessary. The present paper answers these questions within the same analytic class. The contribution is threefold.

The first result, Theorem \ref{thm:obs}, is the observability estimate \eqref{eq:obsL1}, stated with the $L^{1}$-in-time mean square observation norm over $D$ and with a constant that depends on the observation set only through a lower bound $d$ for its measure, and only through an explicit expression in $d^{-1}$. The mean square form \eqref{eq:obs} follows from the Cauchy--Schwarz inequality in time. The constant is explicit in $d$ but is not claimed to be sharp. Compared with \cite{LW26}, the observation set is no longer a cylinder and the equation contains a drift. Compared with the deterministic estimates of \cite{AEW,WZ26}, the adjoint is a backward stochastic equation and the estimate controls the initial state of the adjoint rather than its terminal state.

The second result, Theorem \ref{thm:ctrl}, is null controllability of the forward equation by a single adapted control supported in $D$, with the control bounded in $L^{\infty}$ in time as in \cite{LW26}. This extends the single-control result of \cite{LW26} to a general measurable $D$ and to equations with drift and random initial data.

The third result, Theorem \ref{thm:approx}, is a complete characterization. Within the class \eqref{eq:coef}--\eqref{eq:bound}, equation \eqref{eq:ctrl} is approximately controllable at time $T$ if and only if \eqref{eq:approxcond} holds. For the cylinder $D=G\times(0,T)$ with $\abs{G}>0$ the condition is automatic, which is why this property comes for free in \cite{LW26}. For a general measurable space--time set it does not. In particular, for $D=\T^{n}\times(0,t_{1})$ with $0<t_{1}<T$, Theorem \ref{thm:ctrl} applies while approximate controllability fails.

The reason for this failure is structural. Approximate controllability is the density of the reachable set, and Proposition \ref{prop:orth} identifies the orthogonal complement of that set as the set of terminal data whose adjoint state vanishes at almost every observation time. Density therefore depends on how the observation set accumulates at the terminal time, and \eqref{eq:approxcond} is exactly this accumulation. For a product set $E\times G$ with $G$ open, the same condition reduces to the terminal accumulation of $E$ and was identified in \cite{Lu2011}. The contribution here is the extension to space-dependent analytic coefficients, to a drift coefficient and to an arbitrary measurable space--time set, where the relevant quantity is the measure of $D$ near the terminal time rather than a product structure. The obstruction is the filtration structure of the target, not the regularity of the coefficients or the magnitude of the measure of $D$.

We borrow some ideas from \cite{LW26}, namely analytic energy estimate, propagation of smallness, one-time interpolation and telescoping. The analytic energy estimate, Lemma \ref{lem:scale}, is the counterpart of \cite[Theorem~2.1]{LW26} for equations with a drift and carries the extra factor $\exp(C_{1}T^{1/2})$. The propagation of smallness, Lemma \ref{lem:hilb}, is the Hilbert-valued version of \cite[Theorem~4]{AEW}. The one-time interpolation, Corollary \ref{cor:onet}, is uniform over all measurable sets carrying a prescribed fraction of the measure of the torus. The new ingredient is the reduction of the space--time set $D$ to a measurable set of times on which the spatial slices carry a fixed positive fraction of the torus, followed by the level set subdivision of Lemma \ref{lem:telescope}, which is the criterion of \cite[Proposition~2.2]{WZ26} in the slightly more general form needed here. The extension is needed because the exponent furnished by the one-time interpolation is strictly smaller than one.

We work on the torus, where complex translations are compatible with the Fourier projections and where the energy identity has no boundary term. The analyticity of the spatial profiles is used twice, through the multiplier identity for the translated coefficients and through propagation of smallness. Relaxing it to Gevrey regularity of order larger than one is not covered by this method. 

\subsection*{Organisation}

The rest of the paper is organised as follows. Section \ref{sec:formulation} contains the notation, the main results and the structural results on the reachable set. Section \ref{sec:aux} contains the auxiliary estimates. Section \ref{sec:proof} contains the proofs.

\section{Problem formulation and main results}\label{sec:formulation}

\subsection{Notation and the equations}\label{subsec:notation}

Let $n\ge1$ and $\T^{n}=(\R/2\pi\mathbb{Z})^{n}$. We write $\abs{A}$ for the Lebesgue measure of a measurable set $A$, and $\langle f,g\rangle=\int_{\T^{n}}f\bar g\dd x$ with $\norm{\cdot}$ the associated norm on $L^{2}(\T^{n})$.

Let $(\Omega,\mathcal{F},\mathbb{F},\mathbb{P})$ be a complete filtered probability space with $\mathbb{F}=(\mathcal{F}_{t})_{t\ge0}$ satisfying the usual conditions, carrying a standard scalar Brownian motion $W$, and assume that $\mathcal{F}_{t}=\mathcal{F}_{0}\vee\sigma(W_{s},\,0\le s\le t)$ for every $t\ge0$, with $\mathcal{F}_{0}$ independent of $W$. For a separable Hilbert space $H$ we denote by $L^{2}_{\mathcal{F}_{t}}(\Omega;H)$ the space of $H$-valued $\mathcal{F}_{t}$-measurable random variables $\xi$ with $\E\norm{\xi}^{2}_{H}<\infty$, by $L^{2}_{\mathbb{F}}(0,T;H)$ the space of $H$-valued progressively measurable processes $X$ with $\E\norm{X}^{2}_{L^{2}(0,T;H)}<\infty$, and by $L^{2}_{\mathbb{F}}(\Omega;C([0,T];H))$ the space of $H$-valued adapted processes with continuous paths and $\E\norm{X}^{2}_{C([0,T];H)}<\infty$. For a separable Hilbert space $H_{1}$ and $r\in\{1,\infty\}$ we write $L^{r}_{\mathbb{F}}(0,T;H_{1})$ for the space of $H_{1}$-valued progressively measurable processes $X$ with $\norm{X}_{L^{r}(0,T;H_{1})}<\infty$, the time norm being taken after the $H_{1}$ norm. Thus $u\in L^{\infty}_{\mathbb{F}}(0,T;L^{2}(\Omega;L^{2}(\T^{n})))$ means that $\operatorname{ess\,sup}_{0<t<T}\norm{u(t)}_{L^{2}(\Omega;L^{2}(\T^{n}))}<\infty$. Throughout, $\E$ denotes the expectation with respect to $\mathbb{P}$.

For $\varrho>0$ set $S^{n}_{\varrho}=\{\zeta\in\C^{n}\ \text{such that}\ \abs{\operatorname{Im}\zeta_{j}}<\varrho\ \text{for}\ j=1,\dots,n\}$.
The coefficients are drawn from the following class of analytic functions
\begin{definition}[Analytic class]\label{def:ar}
Let $\varrho>0$. A real-valued function $f$ on $\T^{n}$ belongs to $A_{\varrho}(\T^{n})$ if there is a function $\tilde f$ holomorphic on $S^{n}_{\varrho}$, $2\pi$-periodic in each variable, whose restriction to $\T^{n}$ equals $f$ and which satisfies
\[
\norm{f}_{A_{\varrho}(\T^{n})}:=\sup_{\zeta\in S^{n}_{\varrho}}\abs{\tilde f(\zeta)}<+\infty .
\]
\end{definition}

Let $R>0$ and $C_{1}>0$. The coefficients $a$ and $b$ are real-valued and progressively measurable on $\Omega\times(0,T)\times\T^{n}$. We assume that for almost every $(\omega,t)$ the functions $x\mapsto a(\omega,t,x)$ and $x\mapsto b(\omega,t,x)$ belong to $A_{R}(\T^{n})$, that the functions
\begin{equation}\label{eq:coef}
\alpha(\omega,t)=\norm{a(\omega,t,\cdot)}_{A_{R}(\T^{n})},\qquad \beta(\omega,t)=\norm{b(\omega,t,\cdot)}_{A_{R}(\T^{n})}
\end{equation}
are progressively measurable, and that
\begin{equation}\label{eq:bound}
\norm{\alpha}_{L^{\infty}(\Omega;L^{2}(0,T))}+\norm{\beta}_{L^{\infty}(\Omega;L^{2}(0,T))}\le C_{1}.
\end{equation}
Since $S^{n}_{R}$ contains $\T^{n}$, \eqref{eq:coef} gives $\abs{a(\omega,t,x)}\le\alpha(\omega,t)$ and $\abs{b(\omega,t,x)}\le\beta(\omega,t)$ for every $x\in\T^{n}$.

The first of the two equations is the adjoint, and it reads
\begin{equation}\label{eq:bwd}
\begin{cases}
\dd z+\Delta z\dd t=(a\,z-b\,Z)\dd t+Z\dd W(t), & \text{in }\T^{n}\times(0,T),\\
z(T)=z_{T}\in L^{2}_{\mathcal{F}_{T}}(\Omega;L^{2}(\T^{n})).
\end{cases}
\end{equation}
Under \eqref{eq:coef}--\eqref{eq:bound} the system \eqref{eq:bwd} admits a unique solution $(z,Z)$ in $$\big(L^{2}_{\mathbb{F}}(\Omega;C([0,T];L^{2}(\T^{n})))\cap L^{2}_{\mathbb{F}}(0,T;H^{1}(\T^{n}))\big)\times L^{2}_{\mathbb{F}}(0,T;L^{2}(\T^{n})),$$ 
see \cite[Section~4.4]{LZ}.

The second is the equation to be controlled, and it reads
\begin{equation}\label{eq:fwd}
\begin{cases}
\dd y-\Delta y\dd t=-a\,y\dd t+b\,y\dd W(t), & \text{in }\T^{n}\times(0,T),\\
y(0)=y_{0}\in L^{2}_{\mathcal{F}_{0}}(\Omega;L^{2}(\T^{n})).
\end{cases}
\end{equation}
The control region is an arbitrary measurable subset of the cylinder and is not assumed to be a product. Let $D\subset\T^{n}\times(0,T)$ be measurable and write $D_{t}=\{x\in\T^{n}\ \text{such that}\ (x,t)\in D\}$. For a control $u$ we consider
\begin{equation}\label{eq:ctrl}
\begin{cases}
\dd y-\Delta y\dd t=(-a\,y+\chi_{D}u)\dd t+b\,y\dd W(t), & \text{in }\T^{n}\times(0,T),\\
y(0)=y_{0}\in L^{2}_{\mathcal{F}_{0}}(\Omega;L^{2}(\T^{n})),
\end{cases}
\end{equation}
where $\chi_{D}$ is the characteristic function of $D$ and $u$ is required to vanish outside $D$. For $u\in L^{\infty}_{\mathbb{F}}(0,T;L^{2}(\Omega;L^{2}(\T^{n})))$ the system \eqref{eq:ctrl} admits a unique solution $y\in L^{2}_{\mathbb{F}}(\Omega;C([0,T];L^{2}(\T^{n})))\cap L^{2}_{\mathbb{F}}(0,T;H^{1}(\T^{n}))$ \cite[Theorems~3.18 and 3.20]{LZ}, and the same holds for $u\in L^{2}_{\mathbb{F}}(0,T;L^{2}(\T^{n}))$. The system \eqref{eq:ctrl} is null controllable at time $T$ if for every $y_{0}$ one can choose $u$ so that the corresponding solution satisfies $y(T)=0$ in $L^{2}(\T^{n})$ almost surely, and approximately controllable at time $T$ if for every $y_{0}$, every $y_{1}\in L^{2}_{\mathcal{F}_{T}}(\Omega;L^{2}(\T^{n}))$ and every $\varepsilon>0$ one can choose $u\in L^{2}_{\mathbb{F}}(0,T;L^{2}(\T^{n}))$ vanishing outside $D$ so that $\norm{y(T)-y_{1}}_{L^{2}_{\mathcal{F}_{T}}(\Omega;L^{2}(\T^{n}))}<\varepsilon$.

Let
\begin{equation}\label{eq:weight}
\mu=2C_{1}T^{1/2}+C_{1}^{2},\qquad
h(t)=\sup_{0\le r\le t}e^{\mu r/(2T)}\big(\E\norm{z(r)}^{2}\big)^{1/2},\qquad
\Lambda_{s,t}=1+\abs{t-s}^{-1},
\end{equation}
for a solution $z$ of \eqref{eq:bwd} and $s,t\in[0,T]$ with $s\neq t$.

\subsection{Main results}\label{subsec:results}

Throughout this subsection $a$ and $b$ satisfy \eqref{eq:coef}--\eqref{eq:bound} with the constants $R$ and $C_{1}$. We state the observability estimate first, since the control statements are deduced from it.

\begin{theorem}[Observability from a measurable space--time set]\label{thm:obs}
Let $D\subset\T^{n}\times(0,T)$ be measurable with $\abs{D}\ge d>0$. Then there are $C=C(n,R,C_{1},T)$ and $K=K(n,R)$ such that every solution of \eqref{eq:bwd} satisfies
\begin{equation}\label{eq:obsL1}
\norm{z(0)}_{L^{2}_{\mathcal{F}_{0}}(\Omega;L^{2}(\T^{n}))}\le\exp\big(C(1+d^{-1})^{K}\big)\int_{0}^{T}\norm{z(t)}_{L^{2}_{\mathcal{F}_{t}}(\Omega;L^{2}(D_{t}))}\dd t ,
\end{equation}
and consequently
\begin{equation}\label{eq:obs}
\E\norm{z(0)}^{2}_{L^{2}(\T^{n})}\le\exp\big(C(1+d^{-1})^{K}\big)\E\iint_{D}\abs{z(x,t)}^{2}\dd x\dd t .
\end{equation}
\end{theorem}

The first control statement follows from \eqref{eq:obsL1} by duality, and the constant of the estimate becomes the constant of the control cost.
\begin{theorem}[Null controllability of the forward equation]\label{thm:ctrl}
Let $D\subset\T^{n}\times(0,T)$ be measurable with $\abs{D}\ge d>0$. Then for every $y_{0}\in L^{2}_{\mathcal{F}_{0}}(\Omega;L^{2}(\T^{n}))$ there is a control $u\in L^{\infty}_{\mathbb{F}}(0,T;L^{2}(\Omega;L^{2}(\T^{n})))$ vanishing outside $D$ such that the solution of \eqref{eq:ctrl} satisfies $y(T)=0$ in $L^{2}(\T^{n})$ almost surely and such that
\begin{equation}\label{eq:ctrlcostLinf}
\norm{u}_{L^{\infty}_{\mathbb{F}}(0,T;L^{2}(\Omega;L^{2}(\T^{n})))}\le\exp\big(C(1+d^{-1})^{K}\big)\norm{y_{0}}_{L^{2}_{\mathcal{F}_{0}}(\Omega;L^{2}(\T^{n}))}
\end{equation}
and consequently
\begin{equation}\label{eq:ctrlcost}
\E\norm{u}^{2}_{L^{2}(0,T;L^{2}(\T^{n}))}\le\exp\big(C(1+d^{-1})^{K}\big)\E\norm{y_{0}}^{2}
\end{equation}
with $C=C(n,R,C_{1},T)$ and $K=K(n,R)$.
\end{theorem}

The third statement is a characterization and runs in the opposite direction.
\begin{theorem}[Approximate controllability of the forward equation]\label{thm:approx}
Let $D\subset\T^{n}\times(0,T)$ be measurable. Then \eqref{eq:ctrl} is approximately controllable at time $T$ if and only if
\begin{equation}\label{eq:approxcond}
\abs{D\cap(\T^{n}\times(T-\delta,T))}>0\qquad\text{for every }\delta\in(0,T).
\end{equation}
\end{theorem}

\subsection{The structure of the reachable set}\label{subsec:structure}
The proof of Theorem \ref{thm:approx} rests on the orthogonal complement of the reachable set.
We write
\begin{equation}\label{eq:reach}
\mathcal{Y}=\big\{y_{u}(T)\ \text{with}\ u\in L^{2}_{\mathbb{F}}(0,T;L^{2}(\T^{n}))\ \text{vanishing outside}\ D\big\},
\end{equation}
where $y_{u}$ denotes the solution of \eqref{eq:ctrl} with $y_{0}=0$, and we put
\begin{equation}\label{eq:obstime}
J=\big\{t\in(0,T)\ \text{with}\ \abs{D_{t}}>0\big\}.
\end{equation}
By linearity of \eqref{eq:ctrl}, $\mathcal{Y}$ is a linear subspace of $L^{2}_{\mathcal{F}_{T}}(\Omega;L^{2}(\T^{n}))$, and the terminal values reachable from a given $y_{0}$ form the translate of $\mathcal{Y}$ by the terminal value of the solution of \eqref{eq:ctrl} with $u=0$. Hence \eqref{eq:ctrl} is approximately controllable at time $T$ if and only if $\mathcal{Y}$ is dense in $L^{2}_{\mathcal{F}_{T}}(\Omega;L^{2}(\T^{n}))$.

The next statement identifies this complement in terms of the adjoint state.
\begin{proposition}[Orthogonal complement of the reachable set]\label{prop:orth}
Let $\eta\in L^{2}_{\mathcal{F}_{T}}(\Omega;L^{2}(\T^{n}))$ and let $(z_{\eta},Z_{\eta})$ be the solution of \eqref{eq:bwd} with $z_{T}=\eta$. Let $\mathcal{Y}$ be given by \eqref{eq:reach} and $J$ by \eqref{eq:obstime}. Then $\eta$ is orthogonal to $\mathcal{Y}$ if and only if
\begin{equation}\label{eq:orth}
z_{\eta}(t)=0\quad\text{in }L^{2}_{\mathcal{F}_{t}}(\Omega;L^{2}(\T^{n}))\qquad\text{for almost every }t\in J .
\end{equation}
\end{proposition}

\subsection{Remarks}\label{subsec:remarks}
The remarks concern the coefficient class, the frequency mixing, the direction of time and the scope of the three theorems.
\begin{remark}[The analytic class]\label{rem:class}
Constant coefficients $a\equiv c_{1}$ and $b\equiv c_{2}$ belong to the class \eqref{eq:coef}--\eqref{eq:bound} as soon as $(\abs{c_{1}}+\abs{c_{2}})T^{1/2}\le C_{1}$, so the class is not empty. A Fourier criterion for membership, together with random trigonometric polynomials with arbitrary frequencies, profiles with infinitely many nonzero Fourier modes and profiles whose analytic norm is unbounded in time, is given in \cite{LW26} and is not repeated here. We record three elementary properties of the class and two limitations.

The first property is monotonicity in the width. If $0<\varrho<R$ and $f\in A_{R}(\T^{n})$, then $f\in A_{\varrho}(\T^{n})$ and $\norm{f}_{A_{\varrho}(\T^{n})}\le\norm{f}_{A_{R}(\T^{n})}$, since $S^{n}_{\varrho}\subset S^{n}_{R}$.

The second property is boundedness under translations inside the polystrip. Let $f\in A_{R}(\T^{n})$ with extension $\tilde f$ and let $\varsigma\in\R^{n}$ satisfy $\abs{\varsigma_{j}}<R$ for $j=1,\dots,n$. If $\zeta\in\C^{n}$ satisfies $\abs{\Im\zeta_{j}}+\abs{\varsigma_{j}}<R$ for every $j$, then $\zeta-i\varsigma\in S^{n}_{R}$ and
\[
\abs{\tilde f(\zeta-i\varsigma)}\le\norm{f}_{A_{R}(\T^{n})},
\]
by the definition of $\norm{f}_{A_{R}(\T^{n})}$. Consequently the map $\zeta\mapsto\tilde f(\zeta-i\varsigma)$ is holomorphic and $2\pi$-periodic on the set where $\abs{\Im\zeta_{j}}<R-\abs{\varsigma_{j}}$ for every $j$, and its supremum there is at most $\norm{f}_{A_{R}(\T^{n})}$. The analytic energy estimate uses this property as follows. Multiplication by $b$ is converted into multiplication by a translate of the extension $\tilde b$, and the bound above supplies $\beta(\omega,r)$ as an estimate for every such translate. The frequency mixing is thereby not removed. It is displaced from the spectral projectors to a family of translated multipliers whose size is still measured by $\beta$, and it is this displacement that the energy estimate exploits.

The third property is closure under products. If $f,g\in A_{R}(\T^{n})$, then $fg\in A_{R}(\T^{n})$ and $\norm{fg}_{A_{R}(\T^{n})}\le\norm{f}_{A_{R}(\T^{n})}\norm{g}_{A_{R}(\T^{n})}$, the product of the extensions being holomorphic and $2\pi$-periodic on $S^{n}_{R}$. Hence \eqref{eq:coef}--\eqref{eq:bound} is stable under adding a drift of the same class.

Two limitations are worth recording. The class is uniform, since its elements share a common width and a common bound for the square integral in time of their analytic norms. For instance a profile whose analytic norm behaves like $t^{-1/2}$ near the initial time is analytic at every time and nevertheless outside the class, so \eqref{eq:coef}--\eqref{eq:bound} is strictly stronger than analyticity alone. Smooth profiles which are not analytic are excluded as well, so general bounded space--time potentials as in \cite{WZ26} need not belong to the class.
\end{remark}

\begin{remark}[Frequency mixing]\label{rem:mixing}
Frequency mixing occurs already for the simplest profiles of the class. On $\T^{1}$ take $b(x)=2\cos x$ and let the martingale integrand carry the single spatial mode $x\mapsto e^{i(N+1)x}$ with $N\ge1$. Since
\[
2\cos x\,e^{i(N+1)x}=e^{i(N+2)x}+e^{iNx},
\]
the product of $b$ with that mode carries the mode $e^{iNx}$, whose index lies below the index of the mode carried by the integrand. The low-frequency equation for the adjoint state therefore receives a term generated by the high-frequency part of the martingale integrand, the projected equations do not close, and the commutation with the spectral projectors that underlies the scalar-coefficient argument of \cite{Lu2011,YZ,LWYZ2026} is unavailable.
\end{remark}

\begin{remark}\label{rem:sign}
Theorem \ref{thm:obs} concerns the backward equation. The corresponding observability statement for \eqref{eq:fwd} is not obtained by the method used here. The one-time interpolation is combined with the one-sided energy propagation of Lemma \ref{lem:energy}, which compares the state at an earlier time with the state at a later time in the direction compatible with the backward equation. For \eqref{eq:fwd} the analogue of Lemma \ref{lem:energy} fails. When $a=b=0$ and $y_{0}=e_{k}$ is a Fourier mode, the solution is $y(t)=e^{-t\abs{k}^{2}}e_{k}$, so that the ratio $\E\norm{y(s)}^{2}/\E\norm{y(t)}^{2}$ equals $e^{2\abs{k}^{2}(t-s)}$ for $0\le s<t$, which is unbounded over $k$. This shows that the telescoping argument based on forward energy propagation is not available. It does not rule out observability of \eqref{eq:fwd} by other arguments, a question that is not addressed here. Consequently the forward equation enters this paper only as the object of the control problems in Theorems \ref{thm:ctrl} and \ref{thm:approx}.
\end{remark}

\begin{remark}\label{rem:lw}
Theorems \ref{thm:obs} and \ref{thm:ctrl} differ from \cite[Theorems~1.3 and 1.1]{LW26} in several respects. There the observation region is a cylinder $G\times(0,T)$ over a measurable $G\subset\T^{n}$, the observation norm is the $L^{1}$-in-time mean square norm over that cylinder, and the control space is $L^{\infty}_{\mathbb{F}}(0,T;L^{2}(\Omega;L^{2}(G)))$. Theorems \ref{thm:obs} and \ref{thm:ctrl} keep the $L^{1}$-in-time observation norm and the $L^{\infty}$ control space but allow an arbitrary measurable subset $D$ of the space--time cylinder, together with a drift coefficient and random initial data. In the proof of Theorem \ref{thm:obs}, $D$ is reduced to a measurable set of times on which the spatial slices have measure bounded below by a fixed multiple of $d$, and the one-time interpolation is then applied slice by slice. Theorem \ref{thm:approx} is the counterpart of \cite[Theorem~1.2]{LW26}. For a cylinder the condition \eqref{eq:approxcond} is automatic, so the converse is void, while for a general measurable set the characterization is new. For a product set $E\times G$ with $G$ open and a scalar noise coefficient, the same characterization was obtained in \cite{Lu2011}.
\end{remark}

\begin{remark}\label{rem:dual}
Theorem \ref{thm:ctrl} is obtained from Theorem \ref{thm:obs} by duality, the controlled equation \eqref{eq:ctrl} being paired with \eqref{eq:bwd}. The counterpart for \eqref{eq:bwd}, namely null controllability by a control acting on the drift term, is not established here. Duality would require an observability estimate for \eqref{eq:fwd}, and Remark \ref{rem:sign} explains why the present method does not supply it. Null controllability of \eqref{eq:bwd} is known by other methods for linear and semilinear backward equations, see \cite[Theorem~2.2]{TZ} and \cite[Theorem~1.2]{HLP}. For forward equations, the Carleman-based results of \cite{TZ,HLP} use two controls, one in the drift and one in the diffusion, and hence avoid observability from a single observation region.
\end{remark}

\begin{remark}\label{rem:approxcond}
Condition \eqref{eq:approxcond} depends on $D$ only through the set $J$ of observation times and only through the dichotomy $\abs{D_{t}}>0$ against $\abs{D_{t}}=0$. It is equivalent to $\abs{J\cap(T-\delta,T)}>0$ for every $\delta\in(0,T)$, it is invariant under modification of $D$ on a null set, and it involves neither the spatial shape of the slices nor the magnitude of $\abs{D}$.

The necessary half of Theorem \ref{thm:approx} uses neither the analyticity of $a$ and $b$ nor a sign condition on $b$, but only the bound \eqref{eq:bound}. Assume that \eqref{eq:approxcond} fails, so that $\abs{D\cap(\T^{n}\times(T-\delta,T))}=0$ for some $\delta\in(0,T)$. Take as martingale integrand the indicator of $\T^{n}\times(T-\delta,T)$, let the adjoint state start from zero at time $T-\delta$, and extend both by zero to $(0,T-\delta)$. The resulting pair solves \eqref{eq:bwd} and its terminal value $\eta$ is nonzero, since otherwise uniqueness for \eqref{eq:bwd} would force the pair to vanish while the martingale integrand does not. Moreover the adjoint state vanishes on $D$ outside a set of zero measure. Hence $\eta$ is orthogonal to $\mathcal{Y}$ and $\mathcal{Y}$ is not dense.

The sufficient half uses the analytic class through Proposition \ref{prop:orth}, in passing from a slice of positive measure to the whole torus at each observation time where the slice has positive measure. That passage is the one-time interpolation of Corollary \ref{cor:onet}, applied with the slice as observation set. If $\abs{D_{t}}=\abs{\T^{n}}$ for a sequence $t_{j}\to T$, the passage is void and only the continuity of $t\mapsto z_{\eta}(t)$ at $T$ is needed. Condition \eqref{eq:approxcond} cannot be replaced by the requirement that $T$ be a limit point of $J$, since a countable $J$ gives $\abs{D}=0$ and hence $\mathcal Y=\{0\}$.
\end{remark}

\begin{remark}\label{rem:approxsharp}
Theorem \ref{thm:approx} shows that condition \eqref{eq:approxcond} is indispensable. For $D=\T^{n}\times(0,t_{1})$ with $0<t_{1}<T$ one has $\abs{D}>0$, so that Theorem \ref{thm:ctrl} applies, while $\mathcal{Y}$ is not dense by Theorem \ref{thm:approx}. Hence null controllability at time $T$ does not imply approximate controllability at time $T$ for \eqref{eq:ctrl}. This reproduces the phenomenon identified in \cite{Lu2011} for a scalar noise coefficient, here within the analytic class and with a drift. In particular, the observability estimate of Theorem \ref{thm:obs} does not imply approximate controllability at time $T$. The estimate controls $z(0)$ rather than the terminal datum $z_{T}$, and $z(0)=0$ does not force $z_{T}=0$ for a backward stochastic equation. When $a=b=0$, $z_{T}=W(T)$ and the constant integrand equal to $1$ correspond to the adjoint state $z(t)=W(t)$, which vanishes at $t=0$ and not at $t=T$. Condition \eqref{eq:approxcond} holds in particular when $E\times(T-\delta,T)\subset D$ for some $\delta\in(0,T)$ and some measurable $E\subset\T^{n}$ with $\abs{E}>0$, and it holds for the cylinder $D=G\times(0,T)$ with $\abs{G}>0$. Since bounded adapted controls supported in $G$ are dense in the subspace of $L^{2}_{\mathbb{F}}(0,T;L^{2}(\Omega;L^{2}(\T^{n})))$ consisting of controls vanishing outside $G$, and the map $u\mapsto y_{u}(T)$ is continuous, Theorem \ref{thm:approx} contains the approximate controllability statement of \cite[Theorem~1.2]{LW26} in the case $a=0$ with deterministic $y_{0}$. In the deterministic control problem for the heat equation, with real analytic coefficients, deterministic initial datum and deterministic targets, the reachable set is dense for every measurable $D$ with $\abs{D}>0$, since $\chi_{D}z=0$ forces $z\equiv0$ by unique continuation for real analytic functions. With merely bounded space--time potentials, the same conclusion follows from the observability estimate in \cite{WZ26}. The obstruction isolated here is of a different nature, namely the filtration structure of the target. Even when $a=b=0$ a control supported in $\T^{n}\times(0,t_{1})$ produces a terminal state which is measurable with respect to the noise up to time $t_{1}$, and such states do not exhaust $L^{2}_{\mathcal{F}_{T}}(\Omega;L^{2}(\T^{n}))$. The analytic class enters only in the sufficiency, in the passage from a spatial slice of positive measure to the whole torus at an observation time.
\end{remark}

\section{Auxiliary results}\label{sec:aux}

Throughout this section $a$ and $b$ satisfy \eqref{eq:coef}--\eqref{eq:bound}, and $z$ denotes a solution of \eqref{eq:bwd}.

\subsection{Energy propagation}\label{subsec:energy}

\begin{lemma}[One-sided energy propagation]\label{lem:energy}
Assume \eqref{eq:coef}--\eqref{eq:bound} and let $z$ be the solution of \eqref{eq:bwd}. Let $\mu$ and
$h$ be given by \eqref{eq:weight} with this $z$. Then the following assertions hold.
\begin{enumerate}
\item For all $0\le t_{1}<t_{2}\le T$,
\begin{equation}\label{eq:oneprop}
\E\norm{z(t_{1})}^{2}_{L^{2}(\T^{n})}\le e^{\mu}\E\norm{z(t_{2})}^{2}_{L^{2}(\T^{n})},
\end{equation}
and consequently $\E\norm{z(t)}^{2}_{L^{2}(\T^{n})}\ge e^{-\mu}\E\norm{z(0)}^{2}_{L^{2}(\T^{n})}$ for every
$t\in[0,T]$.
\item The function $h$ is continuous and nondecreasing on $[0,T]$,
$h(0)=\big(\E\norm{z(0)}^{2}_{L^{2}(\T^{n})}\big)^{1/2}$, and
\begin{equation}\label{eq:hcomp}
e^{\mu t/(2T)}\big(\E\norm{z(t)}^{2}_{L^{2}(\T^{n})}\big)^{1/2}
\le h(t)\le e^{\mu/2}e^{\mu t/(2T)}\big(\E\norm{z(t)}^{2}_{L^{2}(\T^{n})}\big)^{1/2}
\end{equation}
for every $t\in[0,T]$. If $\E\norm{z(0)}^{2}_{L^{2}(\T^{n})}>0$, then $h(t)>0$ for every
$t\in[0,T]$.
\end{enumerate}
\end{lemma}

\begin{proof}
For $t\in[0,T]$ write $q(t)=\E\norm{z(t)}^{2}_{L^{2}(\T^{n})}$ and
\[
\Psi(t)=\int_{0}^{t}\big(2\alpha(r)+\beta(r)^{2}\big)\dd r .
\]
The process $\Psi$ is continuous, nondecreasing and $\mathbb{F}$-adapted, and \eqref{eq:bound} together
with the Cauchy--Schwarz inequality gives
\begin{equation}\label{eq:psibound}
\Psi(T)\le2T^{1/2}\norm{\alpha}_{L^{\infty}(\Omega;L^{2}(0,T))}
+\norm{\beta}^{2}_{L^{\infty}(\Omega;L^{2}(0,T))}\le2C_{1}T^{1/2}+C_{1}^{2}=\mu
\end{equation}
almost surely.

Step 1. Let $0\le t_{1}<t_{2}\le T$. Since
$z\in L^{2}_{\mathbb{F}}(\Omega;C([0,T];L^{2}(\T^{n})))\cap L^{2}_{\mathbb{F}}(0,T;H^{1}(\T^{n}))$ and
$Z\in L^{2}_{\mathbb{F}}(0,T;L^{2}(\T^{n}))$, the It\^o formula for $L^{2}(\T^{n})$-valued processes
\cite[Chapter~2]{LZ}, applied to $\norm{z}^{2}$ and \eqref{eq:bwd}, together with
$\langle-\Delta z,z\rangle=\norm{\nabla z}^{2}$ and the identity
$-2\operatorname{Re}\langle Z,bz\rangle+\norm{Z}^{2}=\norm{Z-bz}^{2}-\norm{bz}^{2}$, gives
\begin{equation}\label{eq:energyid}
\norm{z(t_{2})}^{2}-\norm{z(t_{1})}^{2}
=\int_{t_{1}}^{t_{2}}\Big(2\norm{\nabla z}^{2}+2\langle az,z\rangle+\norm{Z-bz}^{2}
-\norm{bz}^{2}\Big)\dd r+\int_{t_{1}}^{t_{2}}2\operatorname{Re}\langle z,Z\rangle\dd W(r).
\end{equation}

Step 2. Fix $\kappa>0$ and let $\tau_{\kappa}$ be the infimum of those $t\in[0,T]$ for which
$\sup_{0\le r\le t}\norm{z(r)}\ge\kappa$, with $\tau_{\kappa}=T$ if there is no such $t$. Since
$\norm{z(r)}\le\kappa$ for $r\le\tau_{\kappa}$,
\[
\E\int_{0}^{T\wedge\tau_{\kappa}}\norm{z}^{2}\norm{Z}^{2}\dd r
\le\kappa^{2}\E\int_{0}^{T}\norm{Z}^{2}\dd r<+\infty ,
\]
so the local martingale in \eqref{eq:energyid}, stopped at $\tau_{\kappa}$, is a martingale. As $\Psi$
is nondecreasing and of finite variation, \eqref{eq:energyid} gives
\[
\begin{aligned}
&e^{\Psi(t_{2}\wedge\tau_{\kappa})}\norm{z(t_{2}\wedge\tau_{\kappa})}^{2}
-e^{\Psi(t_{1}\wedge\tau_{\kappa})}\norm{z(t_{1}\wedge\tau_{\kappa})}^{2}\\
&\quad=\int_{t_{1}\wedge\tau_{\kappa}}^{t_{2}\wedge\tau_{\kappa}}e^{\Psi}
\Big(2\norm{\nabla z}^{2}+\norm{Z-bz}^{2}+2\langle az,z\rangle+2\alpha\norm{z}^{2}
-\norm{bz}^{2}+\beta^{2}\norm{z}^{2}\Big)\dd r+\mathcal M_{\kappa},
\end{aligned}
\]
where $\E\mathcal M_{\kappa}=0$. On the one hand
$\langle az,z\rangle=\int_{\T^{n}}a\abs{z}^{2}\dd x\ge-\alpha\norm{z}^{2}$, and on the other hand
$\norm{bz}^{2}=\int_{\T^{n}}b^{2}\abs{z}^{2}\dd x\le\beta^{2}\norm{z}^{2}$, so the integrand above is
nonnegative. Taking expectations and using $e^{\Psi}\ge1$ and $\Psi\le\mu$ from \eqref{eq:psibound},
\[
\E\norm{z(t_{1}\wedge\tau_{\kappa})}^{2}
\le\E\big[e^{\Psi(t_{1}\wedge\tau_{\kappa})}\norm{z(t_{1}\wedge\tau_{\kappa})}^{2}\big]
\le\E\big[e^{\Psi(t_{2}\wedge\tau_{\kappa})}\norm{z(t_{2}\wedge\tau_{\kappa})}^{2}\big]
\le e^{\mu}\E\norm{z(t_{2}\wedge\tau_{\kappa})}^{2}.
\]
Since $z$ has continuous paths, $\tau_{\kappa}=T$ for every
$\kappa>\sup_{[0,T]}\norm{z(\cdot)}$, and
$\norm{z(t\wedge\tau_{\kappa})}^{2}\le\sup_{[0,T]}\norm{z}^{2}$ with
$\E\sup_{[0,T]}\norm{z}^{2}<+\infty$. Letting $\kappa\to\infty$ by dominated convergence gives
\eqref{eq:oneprop}.

Step 3. The function $q$ is continuous. Hence $h$ is continuous and, by its definition in
\eqref{eq:weight}, nondecreasing, with $h(0)=q(0)^{1/2}$. For $0\le r\le t\le T$, \eqref{eq:oneprop} with
$(t_{1},t_{2})=(r,t)$ gives $e^{\mu r/(2T)}q(r)^{1/2}\le e^{\mu/2}e^{\mu t/(2T)}q(t)^{1/2}$, and taking
the supremum over $r\le t$ yields \eqref{eq:hcomp}. If $q(0)>0$, then \eqref{eq:oneprop} with
$t_{1}=0$ gives $q(t)\ge e^{-\mu}q(0)>0$ for every $t\in[0,T]$, hence $h>0$ on $[0,T]$.
\end{proof}

\subsection{Spatial analyticity}\label{subsec:analyticity}
The energy propagation compares two times but gives no information in the space variable. The following estimate supplies it.
\begin{lemma}[Quantitative spatial analyticity]\label{lem:scale}
Assume \eqref{eq:coef}--\eqref{eq:bound} and let $(z,Z)$ be the solution of \eqref{eq:bwd}. Then for every $0\le t<s\le T$ the random variable $z(t)$ is the trace on $\T^{n}$ of a function $\mathcal Z_{t}$ holomorphic on $S^{n}_{3R/4}$, $2\pi$-periodic in each variable, with values in $L^{2}_{\mathcal{F}_{t}}(\Omega;\C)$, and
\begin{equation}\label{eq:analytic}
\sup_{\zeta\in S^{n}_{R/2}}\norm{\mathcal Z_{t}(\zeta)}_{L^{2}(\Omega)}\le C\exp\Big(\frac{C_{1}^{2}}{2}+C_{1}T^{1/2}+\frac{9nR^{2}}{64(s-t)}\Big)\norm{z(s)}_{L^{2}_{\mathcal{F}_{s}}(\Omega;L^{2}(\T^{n}))},
\end{equation}
where $C=C(n,R)$.
\end{lemma}

\begin{proof}
Let $e_{k}(x)=(2\pi)^{-n/2}e^{ik\cdot x}$ for $k\in\mathbb{Z}^{n}$, let $\abs{k}_{1}=\sum_{j=1}^{n}\abs{k_{j}}$ and $\abs{k}_{\infty}=\max_{1\le j\le n}\abs{k_{j}}$, and put $D_{j}=-i\partial_{x_{j}}$ and $\sigma\cdot D=\sum_{j=1}^{n}\sigma_{j}D_{j}$ for $\sigma\in\{+1,-1\}^{n}$. For $N\in\mathbb{N}$ let $P_{N}$ be the orthogonal projection of $L^{2}(\T^{n})$ onto $H_{N}=\operatorname{span}\{e_{k}\ \text{with}\ \abs{k}_{\infty}\le N\}$, and let $(z_{N},Z_{N})$ be the $H_{N}$-valued solution of
\begin{equation}\label{eq:galerkin}
\dd z_{N}=-(\Delta z_{N}+P_{N}(bZ_{N}))\dd r+P_{N}(az_{N})\dd r+Z_{N}\dd W(r),\qquad z_{N}(s)=P_{N}z(s),
\end{equation}
on $[0,s]$. Set $e_{N}=z_{N}-P_{N}z$ and $h_{N}=Z_{N}-P_{N}Z$. Subtracting from \eqref{eq:galerkin} the projection of \eqref{eq:bwd} onto $H_{N}$ gives
\[
\dd e_{N}=\Big(-\Delta e_{N}+P_{N}(ae_{N})-P_{N}\big(a(I-P_{N})z\big)-P_{N}(bh_{N})+P_{N}\big(b(I-P_{N})Z\big)\Big)\dd r+h_{N}\dd W(r),\ e_{N}(s)=0 .
\]
Since $\abs{a}\le\alpha$ and $\abs{b}\le\beta$, the proof of \cite[Appendix A]{LW26} carries over to this equation once the weight built there from $\beta$ is replaced by $\Gamma_{0}(r)=\exp\big(\int_{0}^{r}(2\alpha+2\beta^{2})\dd\tau\big)$, which is bounded by $\exp(2C_{1}T^{1/2}+2C_{1}^{2})$ almost surely by \eqref{eq:bound}. The two forcing terms in the previous display are bounded by $\alpha\norm{(I-P_{N})z}+\beta\norm{(I-P_{N})Z}$, and by the Cauchy--Schwarz inequality in time, \eqref{eq:bound} and the strong convergence of $P_{N}$ they satisfy
\[
\E\Big(\int_{0}^{s}\big(\alpha\norm{(I-P_{N})z}+\beta\norm{(I-P_{N})Z}\big)\dd r\Big)^{2}\longrightarrow0
\]
as $N\to\infty$. Hence $z_{N}\to z$ in $L^{2}_{\mathbb{F}}(\Omega;C([0,s];L^{2}(\T^{n})))\cap L^{2}_{\mathbb{F}}(0,s;H^{1}(\T^{n}))$ and $Z_{N}\to Z$ in $L^{2}_{\mathbb{F}}(0,s;L^{2}(\T^{n}))$.

Step 1. Fix $\varrho\in(0,R)$ and a sign vector $\sigma\in\{+1,-1\}^{n}$. Put $\delta=s-t$ and $\theta(r)=\varrho(s-r)/\delta$ for $r\in[t,s]$, so that $\theta(t)=\varrho$, $\theta(s)=0$ and $\theta'=-\varrho/\delta$. Let $\tilde a$ and $\tilde b$ denote the holomorphic extensions of $a$ and $b$ to $S^{n}_{R}$ furnished by Definition \ref{def:ar}. These extensions are unique by analytic continuation, and we choose the representatives that are measurable with respect to the product of the progressive $\sigma$-algebra on $\Omega\times(0,T)$ and the Borel $\sigma$-algebra of $S^{n}_{R}$. Such a choice is possible since on $S^{n}_{R}$ the extensions coincide with the Fourier series in $x$ of $a$ and $b$. For $v\in H_{N}$ the multiplier $e^{\theta\sigma\cdot D}$ acts by $e^{\theta\sigma\cdot D}e_{k}=e^{\theta\sigma\cdot k}e_{k}$ and satisfies $(e^{\theta\sigma\cdot D}v)(x)=v(x-i\theta\sigma)$. Since $P_{N}$ commutes with Fourier multipliers, setting $b_{\sigma}(r,x)=\tilde b(r,x-i\theta(r)\sigma)$ and $a_{\sigma}(r,x)=\tilde a(r,x-i\theta(r)\sigma)$ gives, by Cauchy's theorem in each coordinate as in \cite[eq.~(2.6)]{LW26}, for $v,w\in H_{N}$,
\[
e^{\theta\sigma\cdot D}P_{N}(bw)=P_{N}\big(b_{\sigma}e^{\theta\sigma\cdot D}w\big),\qquad e^{\theta\sigma\cdot D}P_{N}(av)=P_{N}\big(a_{\sigma}e^{\theta\sigma\cdot D}v\big),
\]
and $a_{\sigma}$ and $b_{\sigma}$ are progressively measurable. Since $\abs{\theta(r)}\le\varrho<R$, one has $\abs{b_{\sigma}}\le\beta$ and $\abs{a_{\sigma}}\le\alpha$. Set $X_{N}=e^{\theta\sigma\cdot D}z_{N}$ and $Y_{N}=e^{\theta\sigma\cdot D}Z_{N}$. Applying \eqref{eq:galerkin} and It\^o's formula,
\[
\dd\norm{X_{N}}^{2}=\Big(2\norm{\nabla X_{N}}^{2}+2\theta'\operatorname{Re}\langle\sigma\cdot DX_{N},X_{N}\rangle+\norm{Y_{N}}^{2}-2\operatorname{Re}\langle X_{N},b_{\sigma}Y_{N}\rangle+2\operatorname{Re}\langle X_{N},a_{\sigma}X_{N}\rangle\Big)\dd r+\dd M_{N}
\]
with $M_{N}$ a martingale. Completing the squares,
\[
2\norm{\nabla X_{N}}^{2}+2\theta'\operatorname{Re}\langle\sigma\cdot DX_{N},X_{N}\rangle=2\sum_{j=1}^{n}\Big\|\Big(D_{j}+\frac{\theta'\sigma_{j}}{2}\Big)X_{N}\Big\|^{2}-\frac{n\abs{\theta'}^{2}}{2}\norm{X_{N}}^{2},
\]
\[
\norm{Y_{N}}^{2}-2\operatorname{Re}\langle X_{N},b_{\sigma}Y_{N}\rangle=\norm{Y_{N}-\overline{b_{\sigma}}X_{N}}^{2}-\norm{b_{\sigma}X_{N}}^{2}.
\]
Let $\Gamma(r)=\exp\big(\int_{t}^{r}\big(n\abs{\theta'}^{2}/2+2\alpha+\beta^{2}\big)\dd\tau\big)$. Since $\norm{b_{\sigma}X_{N}}^{2}\le\beta^{2}\norm{X_{N}}^{2}$ and $2\operatorname{Re}\langle X_{N},a_{\sigma}X_{N}\rangle\ge-2\alpha\norm{X_{N}}^{2}$, the drift of $\Gamma\norm{X_{N}}^{2}$ is nonnegative. The stochastic integral has zero expectation by the integrability of $z_{N}$ and $Z_{N}$ on the finite dimensional space $H_{N}$, see \cite[Theorem~2.130]{LZ}. Integrating over $[t,s]$ and using $\Gamma(t)=1$, $X_{N}(s)=P_{N}z(s)$, $\int_{t}^{s}n\abs{\theta'}^{2}/2=n\varrho^{2}/(2(s-t))$, $\int_{t}^{s}\beta^{2}\le C_{1}^{2}$ and $\int_{t}^{s}2\alpha\le2C_{1}T^{1/2}$ gives
\[
\E\norm{e^{\varrho\sigma\cdot D}z_{N}(t)}^{2}\le\exp\Big(\frac{n\varrho^{2}}{2(s-t)}+C_{1}^{2}+2C_{1}T^{1/2}\Big)\E\norm{z(s)}^{2}.
\]

Step 2. For each $k\in\mathbb{Z}^{n}$ choose $\sigma$ with $\sigma\cdot k=\abs{k}_{1}$. Then $e^{2\varrho\abs{k}_{1}}\le\sum_{\sigma}e^{2\varrho\sigma\cdot k}$, and Parseval's identity gives, for the multiplier $e^{\varrho\abs{D}_{1}}v=\sum_{k}e^{\varrho\abs{k}_{1}}\hat v(k)e_{k}$,
\[
\E\norm{e^{\varrho\abs{D}_{1}}z_{N}(t)}^{2}\le2^{n}\exp\Big(\frac{n\varrho^{2}}{2(s-t)}+C_{1}^{2}+2C_{1}T^{1/2}\Big)\E\norm{z(s)}^{2}.
\]
Since $P_{K}z_{N}(t)\to P_{K}z(t)$ in $L^{2}_{\mathcal{F}_{t}}(\Omega;L^{2}(\T^{n}))$ and $e^{\varrho\abs{D}_{1}}$ is bounded on $H_{K}$, letting $N\to\infty$ and then $K\to\infty$ with monotone convergence yields
\begin{equation}\label{eq:expfreq}
\E\norm{e^{\varrho\abs{D}_{1}}z(t)}^{2}\le2^{n}\exp\Big(\frac{n\varrho^{2}}{2(s-t)}+C_{1}^{2}+2C_{1}T^{1/2}\Big)\E\norm{z(s)}^{2}.
\end{equation}

Step 3. Write $z(t)=\sum_{k}z_{k}(t)e_{k}$ with $z_{k}(t)\in L^{2}_{\mathcal{F}_{t}}(\Omega;\C)$. For $0<\varrho_{0}<\varrho$ and $\zeta\in S^{n}_{\varrho_{0}}$ the Cauchy--Schwarz inequality gives
\[
\sum_{k\in\mathbb{Z}^{n}}\norm{z_{k}(t)e^{ik\cdot\zeta}}_{L^{2}(\Omega)}\le\Big(\sum_{k\in\mathbb{Z}^{n}}e^{-2(\varrho-\varrho_{0})\abs{k}_{1}}\Big)^{1/2}\Big(\sum_{k\in\mathbb{Z}^{n}}e^{2\varrho\abs{k}_{1}}\norm{z_{k}(t)}^{2}_{L^{2}(\Omega)}\Big)^{1/2},
\]
and the right hand side is finite by \eqref{eq:expfreq}. Hence the series converges absolutely and uniformly on $S^{n}_{\varrho_{0}}$ with values in $L^{2}_{\mathcal{F}_{t}}(\Omega;\C)$, its sum $\mathcal Z_{t}$ is holomorphic and periodic on $S^{n}_{\varrho_{0}}$, and $\mathcal Z_{t}=z(t)$ on $\T^{n}$. The same estimate yields
\[
\sup_{\zeta\in S^{n}_{\varrho_{0}}}\norm{\mathcal Z_{t}(\zeta)}_{L^{2}(\Omega)}\le C(n,\varrho-\varrho_{0})\exp\Big(\frac{n\varrho^{2}}{4(s-t)}+\frac{C_{1}^{2}}{2}+C_{1}T^{1/2}\Big)\norm{z(s)}_{L^{2}_{\mathcal{F}_{s}}(\Omega;L^{2}(\T^{n}))}.
\]
Since $0<\varrho_{0}<\varrho<R$ are arbitrary and the series defining $\mathcal Z_{t}$ does not depend on $\varrho$, the function $\mathcal Z_{t}$ is holomorphic on $S^{n}_{R}$, hence on $S^{n}_{3R/4}$. Taking $\varrho=3R/4$ and $\varrho_{0}=R/2$ in the last estimate gives \eqref{eq:analytic}.
\end{proof}

\begin{remark}\label{rem:scalecmp}
Lemma \ref{lem:scale} is the counterpart of \cite[Theorem~2.1]{LW26} for equations with a drift term. The factor $\exp(C_{1}T^{1/2})$ is the cost of the drift, and only the exponential dependence on $(s-t)^{-1}$ in \eqref{eq:analytic} enters the subsequent interpolation and telescoping arguments.
\end{remark}

\subsection{Propagation of smallness and interpolation}\label{subsec:interp}

\begin{lemma}[Hilbert valued analytic propagation]\label{lem:hilb}
Let $\mathcal H$ be a separable complex Hilbert space, let $r>0$ and let $E\subset\T^{n}$ be measurable with $\abs{E}>0$, and put $\nu=\abs{E}/\abs{\T^{n}}$. Then there is $\vartheta\in(0,1)$, depending only on $n$, $r$ and $\nu$, and there are $\kappa_{1}=\kappa_{1}(n,r)$, $\kappa_{2}=\kappa_{2}(n,r)$ and $\kappa_{3}=\kappa_{3}(n,r)$ such that $\vartheta^{-1}\le\kappa_{1}\big(1+\log(1/\nu)\big)$ and such that for every $\mathcal H$-valued function $f$ holomorphic on $S^{n}_{r}$, $2\pi$-periodic in each variable, with $M:=\sup_{\zeta\in S^{n}_{r}}\norm{f(\zeta)}_{\mathcal H}<+\infty$,
\begin{equation}\label{eq:hilb}
\norm{f}_{L^{2}(\T^{n};\mathcal H)}\le\kappa_{2}\nu^{-\kappa_{3}}\norm{f}^{\vartheta}_{L^{2}(E;\mathcal H)}M^{1-\vartheta}.
\end{equation}
\end{lemma}

\begin{proof}
If $M=0$ or $\mathcal H=\{0\}$ the claim is immediate. Assume $M>0$ and $\mathcal H\neq\{0\}$. Regard $E$ as a subset of $Q=[-\pi,\pi)^{n}$ and set $\varrho=\min\{1,r/(4\pi\sqrt{n})\}$. Then $E\subset Q\subset B_{2\pi\sqrt{n}}(0)$ and $2\pi\sqrt{n}\varrho\le r/2$. For a unit vector $\xi\in\mathcal H$ and $\gamma\in\R$ put $\varphi(\zeta)=\langle f(\zeta),\xi\rangle_{\mathcal H}$ for $\zeta\in S^{n}_{r}$ and $g(x)=\operatorname{Re}\big(e^{-i\gamma}\varphi(x)\big)$ for $x\in\R^{n}$. The map $\varphi$ is holomorphic on $S^{n}_{r}$ and bounded by $M$, and every polydisc of radius $r/2$ centred at a real point lies in $S^{n}_{r}$. Cauchy's inequality for $\varphi$ on this polydisc, together with $\partial^{\kappa}g=\operatorname{Re}\partial^{\kappa}\varphi$ at real points, gives, for $x\in B_{4\pi\sqrt{n}}(0)$ and every multi-index $\kappa$,
\[
\abs{\partial^{\kappa}g(x)}\le M\kappa!\Big(\frac{2}{r}\Big)^{\abs{\kappa}}\le M\frac{\abs{\kappa}!}{(2\pi\sqrt{n}\varrho)^{\abs{\kappa}}}.
\]
By \cite[Theorem 4]{AEW} there are $N\ge1$ and $\vartheta\in(0,1)$, depending only on $n$, $\varrho$ and $\abs{E}/\abs{B_{2\pi\sqrt{n}}(0)}$, such that
\[
\sup_{x\in B_{2\pi\sqrt{n}}(0)}\abs{g(x)}\le N\Big(\frac{1}{\abs{E}}\int_{E}\abs{g(x)}\dd x\Big)^{\vartheta}M^{1-\vartheta}.
\]
The dependence of $N$ and $\vartheta$ on $\abs{E}$ is explicit in the proof. The one dimensional step is \cite[(3.1) and Lemma~3.2]{AE2013}. There a function $\psi$ holomorphic in the unit disc with $\abs{\psi}\le1$ and a measurable set $F$ of the real axis satisfy
\[
\sup_{B_{1/2}}\abs{\psi}\le\Big(\sup_{F}\abs{\psi}\Big)\Big(\frac{3}{\abs{F}}\Big)^{j}+2\Big(\frac{7}{8}\Big)^{j}
\]
for every integer $j\ge1$. Put $\eta=\sup_{F}\abs{\psi}$ and
\[
\theta_{F}=\frac{\log(8/7)}{\log(3/\abs{F})+\log(8/7)}.
\]
If $\eta\le7\abs{F}/12$, take $j$ equal to the integer part of
$\log(2/\eta)\big/\big(\log(3/\abs{F})+\log(8/7)\big)$, and take $j=1$ otherwise. Then
\[
\sup_{B_{1/2}}\abs{\psi}\le6\abs{F}^{-1}\eta^{\theta_{F}},
\qquad
\theta_{F}^{-1}=1+\frac{\log(3/\abs{F})}{\log(8/7)}\le9.3\Big(1+\log\frac{1}{\abs{F}}\Big).
\]
Restricting $f$ to the real lines of a fixed direction replaces $F$ by a set of measure at least $c(n)\abs{E}$, and the iteration of Hadamard's three circle theorem in \cite[Lemma~3.3]{AE2013} uses a number of steps depending only on $\varrho$. Hence $\vartheta^{-1}\le\kappa_{1}(1+\log(1/\nu))$ and $N\le\kappa_{2}\nu^{-\kappa_{3}}$, with $\kappa_{1},\kappa_{2},\kappa_{3}$ depending only on $n$ and $r$.
Since $\int_{E}\abs{g}\le\int_{E}\norm{f}_{\mathcal H}\le\abs{E}^{1/2}\norm{f}_{L^{2}(E;\mathcal H)}$, we obtain
\[
\sup_{x\in Q}\abs{g(x)}\le\kappa_{2}\nu^{-\kappa_{3}}\norm{f}^{\vartheta}_{L^{2}(E;\mathcal H)}M^{1-\vartheta},
\]
where the factor $\abs{E}^{-\vartheta/2}$ has been absorbed into $\kappa_{2}$. This is legitimate because $\abs{E}=\nu\abs{\T^{n}}$ with $\abs{\T^{n}}>1$ gives $\log(1/\abs{E})\le\log(1/\nu)$, hence $\vartheta\log(1/\abs{E})\le1/\kappa_{1}$. Taking the supremum over unit vectors $\xi$ and over $\gamma$, and using the dual characterisation of the norm,
\[
\sup_{x\in Q}\norm{f(x)}_{\mathcal H}=\sup_{\norm{\xi}_{\mathcal H}=1}\sup_{\gamma\in\R}\sup_{x\in Q}\abs{g(x)}\le\kappa_{2}\nu^{-\kappa_{3}}\norm{f}^{\vartheta}_{L^{2}(E;\mathcal H)}M^{1-\vartheta}.
\]
Since $f$ is $2\pi$-periodic, integrating over $Q$ gives \eqref{eq:hilb} after absorbing $\abs{Q}^{1/2}$ into $\kappa_{2}$.
\end{proof}

Combining Lemma \ref{lem:scale} with Lemma \ref{lem:hilb} gives an inequality at a single observation time.
\begin{corollary}[One time interpolation]\label{cor:onet}
Assume \eqref{eq:coef}--\eqref{eq:bound} and let $\nu\in(0,1]$. Then there are $\vartheta\in(0,1)$ and $C_{\nu}\ge1$, depending only on $n$, $R$, $C_{1}$, $T$ and $\nu$, such that for every $0\le t<s\le T$, every measurable $E\subset\T^{n}$ with $\abs{E}\ge\nu\abs{\T^{n}}$ and every solution of \eqref{eq:bwd} one has
\begin{equation}\label{eq:onet}
\norm{z(t)}_{L^{2}_{\mathcal{F}_{t}}(\Omega;L^{2}(\T^{n}))}\le C_{\nu}\exp\Big(\frac{C_{\nu}}{s-t}\Big)\norm{z(t)}^{\vartheta}_{L^{2}_{\mathcal{F}_{t}}(\Omega;L^{2}(E))}\norm{z(s)}^{1-\vartheta}_{L^{2}_{\mathcal{F}_{s}}(\Omega;L^{2}(\T^{n}))}.
\end{equation}
\end{corollary}

\begin{proof}
Let $\mathcal H=L^{2}_{\mathcal{F}_{t}}(\Omega;\C)$ and $f=\mathcal Z_{t}$. By Lemma \ref{lem:scale}, $f$ is $\mathcal H$-valued and holomorphic on $S^{n}_{3R/4}$, hence on $S^{n}_{R/2}$, and
\[
M:=\sup_{\zeta\in S^{n}_{R/2}}\norm{f(\zeta)}_{\mathcal H}
\le C(n,R)\exp\Big(\frac{C_{1}^{2}}{2}+C_{1}T^{1/2}+\frac{9nR^{2}}{64(s-t)}\Big)
\norm{z(s)}_{L^{2}_{\mathcal{F}_{s}}(\Omega;L^{2}(\T^{n}))}.
\]
Fix a measurable $E_{0}\subset E$ with $\abs{E_{0}}=\nu\abs{\T^{n}}$, which exists since $r\mapsto\abs{E\cap[-r,r]^{n}}$ is continuous from $0$ to $\abs{E}$. Since $\abs{E_{0}}=\nu\abs{\T^{n}}$, Lemma \ref{lem:hilb} with $r=R/2$ and $E_{0}$ applies with the same parameter $\nu$ as in the statement, and gives $\vartheta\in(0,1)$ with $\vartheta^{-1}\le\kappa_{1}(1+\log(1/\nu))$ and $C_{h}=\kappa_{2}\nu^{-\kappa_{3}}$ such that
\[
\norm{f}_{L^{2}(\T^{n};\mathcal H)}
\le C_{h}\norm{f}^{\vartheta}_{L^{2}(E_{0};\mathcal H)}M^{1-\vartheta}.
\]
By Fubini's theorem,
\[
\norm{f}^{2}_{L^{2}(\T^{n};\mathcal H)}
=\int_{\T^{n}}\E\abs{z(t,x)}^{2}\dd x
=\norm{z(t)}^{2}_{L^{2}_{\mathcal{F}_{t}}(\Omega;L^{2}(\T^{n}))},
\]
and the same identity holds with $E_{0}$ in place of $\T^{n}$. Since $E_{0}\subset E$,
\[
\norm{z(t)}_{L^{2}_{\mathcal{F}_{t}}(\Omega;L^{2}(E_{0}))}
\le\norm{z(t)}_{L^{2}_{\mathcal{F}_{t}}(\Omega;L^{2}(E))}.
\]
Combining the displayed estimates gives
\[
\begin{aligned}
\norm{z(t)}_{L^{2}_{\mathcal{F}_{t}}(\Omega;L^{2}(\T^{n}))}
&\le C_{h}C(n,R)^{1-\vartheta}
\exp\Big((1-\vartheta)\Big(\frac{C_{1}^{2}}{2}+C_{1}T^{1/2}\Big)\Big)\\
&\quad\times\exp\Big(\frac{(1-\vartheta)9nR^{2}}{64(s-t)}\Big)
\norm{z(t)}^{\vartheta}_{L^{2}_{\mathcal{F}_{t}}(\Omega;L^{2}(E))}
\norm{z(s)}^{1-\vartheta}_{L^{2}_{\mathcal{F}_{s}}(\Omega;L^{2}(\T^{n}))}.
\end{aligned}
\]
Set
\[
C_{\nu}:=\max\Big\{1,\;C_{h}C(n,R)^{1-\vartheta}
\exp\Big((1-\vartheta)\Big(\frac{C_{1}^{2}}{2}+C_{1}T^{1/2}\Big)\Big),\;\frac{9nR^{2}}{64}\Big\}.
\]
Then $C_{\nu}$ depends only on $n$, $R$, $C_{1}$, $T$ and $\nu$, and
\[
\exp\Big(\frac{(1-\vartheta)9nR^{2}}{64(s-t)}\Big)
\le\exp\Big(\frac{C_{\nu}}{s-t}\Big).
\]
Therefore the last estimate yields \eqref{eq:onet}.
\end{proof}

\subsection{Telescoping over a measurable time set}\label{subsec:telescope}
It remains to iterate \eqref{eq:onet} over the observation times.
\begin{lemma}[Telescoping over a measurable time set]\label{lem:telescope}
Let $T>0$, $1<\gamma<2$, $0<\delta\le1$, $C_{0}\ge1$ and $a_{0}\in(0,T)$. Let $h\in C([0,T])$ be nonnegative and nonincreasing, let $A\subset(0,T)$ be measurable with $\abs{A}\ge a_{0}$, and let $g\in L^{1}(A)$ be nonnegative. Assume that
\begin{equation}\label{eq:teleshyp}
h(t)\le C_{0}\exp\Big(C_{0}\Big(1+(t-s)^{-1}+\log\frac{h(s)}{h(t)}\Big)^{\gamma}\Big)g(t)^{\delta}h(s)^{1-\delta}
\end{equation}
for all $0\le s<t\le T$ and almost every $t\in A$ such that $h(t)>0$. Then
\begin{equation}\label{eq:teles}
h(T)\le C\int_{A}g(t)\dd t,\qquad C=C(C_{0},\gamma,\delta,T,a_{0}).
\end{equation}
\end{lemma}

\begin{proof}
The argument is that of \cite[Proposition~2.2]{WZ26}, where the case $\delta=1$ is treated. The hypothesis \eqref{eq:teleshyp} is weaker, and the only modification occurs in Step 1 below, where the factor $h(s)^{1-\delta}$ is absorbed by a large multiple of $h(s)$; Steps 2 and 3 are those of \cite[Proposition~2.2]{WZ26} and are reproduced for completeness.

If $h(T)=0$ there is nothing to prove. Assume $h(T)>0$. Since $h$ is nonnegative and nonincreasing, $h(t)\ge h(T)>0$ for every $t\in[0,T]$, and \eqref{eq:teleshyp} holds for almost every $t\in A$.

Step 1. Fix $\eta\in(0,1)$ and put $\lambda(\eta)=C_{0}^{1/\delta}\eta^{-(1-\delta)/\delta}$. We claim that
\begin{equation}\label{eq:teleslin}
h(t)\le\eta h(s)+\lambda(\eta)\exp\Big(\frac{C_{0}}{\delta}\Big(1+(t-s)^{-1}+\log\frac{1}{\eta}\Big)^{\gamma}\Big)g(t)
\end{equation}
for all $0\le s<t\le T$ and almost every $t\in A$. If $h(t)\le\eta h(s)$, then \eqref{eq:teleslin} holds because $g\ge0$. If $h(t)>\eta h(s)$, then $h(s)/h(t)<1/\eta$, and dividing \eqref{eq:teleshyp} by $h(t)^{1-\delta}$ followed by the power $1/\delta$ gives
\[
h(t)\le C_{0}^{1/\delta}\exp\Big(\frac{C_{0}}{\delta}\Big(1+(t-s)^{-1}+\log\frac{h(s)}{h(t)}\Big)^{\gamma}\Big)\Big(\frac{h(s)}{h(t)}\Big)^{(1-\delta)/\delta}g(t).
\]
The last factor is at most $\eta^{-(1-\delta)/\delta}$ and the exponent is at most $\frac{C_{0}}{\delta}(1+(t-s)^{-1}+\log(1/\eta))^{\gamma}$, so \eqref{eq:teleslin} follows from $\eta h(s)\ge0$.

Step 2. Define $\Phi(t)=\abs{A\cap(0,t)}$ for $t\in[0,T]$. Then $\Phi$ is continuous and nondecreasing, $\Phi(0)=0$, $\Phi(T)=\abs{A}$ and $\Phi'=\chi_{A}$ almost everywhere. Choose $p=p(\gamma)$ and $q=q(\gamma)$ with
\[
1<p<\frac{1}{\gamma-1},\qquad p\gamma<q<\frac{\gamma}{\gamma-1},
\]
which is possible since $1<\gamma<2$. For $k\ge1$ put $m_{k}=\frac{\abs{A}}{2}k^{1-p}$, $\tau_{k}=\inf\{t\in[0,T]:\Phi(t)\ge m_{k}\}$ and $\ell_{k}=m_{k}-m_{k+1}$. Since $p>1$ the sequence $m_{k}$ decreases to $0$, so $\tau_{k}$ decreases and $\Phi(\tau_{k})=m_{k}$ with $\tau_{k}\in(0,T)$. Since $x\mapsto x^{1-p}$ is convex and decreasing on $(0,\infty)$,
\[
\ell_{k}=\frac{\abs{A}}{2}\big(k^{1-p}-(k+1)^{1-p}\big)\ge c_{0}k^{-p},\qquad c_{0}=\frac{a_{0}}{2}(p-1)2^{-p},
\]
where $\abs{A}\ge a_{0}$ was used. Put
\[
A_{k}=\big\{t\in A\cap(\tau_{k+1},\tau_{k}):\Phi(t)-m_{k+1}\ge\ell_{k}/2\big\}.
\]
Since $\Phi$ increases at rate $1$ on $A$, one has $\abs{A_{k}}=\ell_{k}/2$ and $t-\tau_{k+1}\ge\ell_{k}/2$ for every $t\in A_{k}$. Applying \eqref{eq:teleslin} with $\eta=\eta_{k}$, $s=\tau_{k+1}$ and $t\in A_{k}$, and using $h(\tau_{k})\le h(t)$ and $(t-\tau_{k+1})^{-1}\le2/\ell_{k}$, gives, after averaging over $A_{k}$,
\begin{equation}\label{eq:telesstep}
h(\tau_{k})\le\eta_{k}h(\tau_{k+1})+\frac{2}{\ell_{k}}\lambda(\eta_{k})\exp\Big(\frac{C_{0}}{\delta}\Big(1+\frac{2}{\ell_{k}}+\log\frac{1}{\eta_{k}}\Big)^{\gamma}\Big)\int_{A_{k}}g(t)\dd t .
\end{equation}

Step 3. Set $w_{k}=e^{-k^{q}}$ and $\eta_{k}=w_{k+1}/w_{k}\in(0,1)$ for $k\ge1$. Then
\[
\log\frac{1}{\eta_{k}}=(k+1)^{q}-k^{q}\le q2^{q-1}k^{q-1},\qquad
\log\lambda(\eta_{k})\le\frac{1}{\delta}\log C_{0}+\frac{1-\delta}{\delta}q2^{q-1}k^{q-1},
\]
and, with a constant $C_{\gamma}$ depending only on $\gamma$,
\[
\Big(1+\frac{2}{\ell_{k}}+\log\frac{1}{\eta_{k}}\Big)^{\gamma}\le C_{\gamma}\big(1+k^{p\gamma}+k^{\gamma(q-1)}\big).
\]
Let $\Theta_{k}$ be the coefficient of $\int_{A_{k}}g$ in \eqref{eq:telesstep} multiplied by $w_{k}$. Using $\ell_{k}^{-1}\le c_{0}^{-1}k^{p}$ and $w_{k}\eta_{k}=w_{k+1}$,
\[
\log\Theta_{k}\le-k^{q}+C_{1}'k^{q-1}+C_{2}'k^{p\gamma}+C_{3}'k^{\gamma(q-1)}+C_{4}'\log k+C_{5}',
\]
where $C_{1}',\dots,C_{5}'$ depend only on $C_{0}$, $\gamma$, $\delta$ and $a_{0}$. Since $p\gamma<q$, $\gamma(q-1)<q$ and $q-1<q$, there is $k_{0}$, depending only on $C_{0}$, $\gamma$, $\delta$ and $a_{0}$, such that $\Theta_{k}\le1$ for every $k\ge k_{0}$.

The sets $A_{k}$ are pairwise disjoint subsets of $A$. Multiplying \eqref{eq:telesstep} by $w_{k}$ and summing over $k_{0}\le k\le K$ therefore yields
\[
w_{k_{0}}h(\tau_{k_{0}})-w_{K+1}h(\tau_{K+1})=\sum_{k=k_{0}}^{K}\big(w_{k}h(\tau_{k})-w_{k+1}h(\tau_{k+1})\big)\le\sum_{k=k_{0}}^{K}\int_{A_{k}}g(t)\dd t\le\int_{A}g(t)\dd t .
\]
Since $w_{K+1}\to0$ and $0\le h\le h(0)<+\infty$, letting $K\to\infty$ gives $w_{k_{0}}h(\tau_{k_{0}})\le\int_{A}g$. Finally $h(T)\le h(\tau_{k_{0}})$ because $\tau_{k_{0}}\le T$ and $h$ is nonincreasing, hence
\[
h(T)\le w_{k_{0}}^{-1}\int_{A}g(t)\dd t=e^{k_{0}^{q}}\int_{A}g(t)\dd t ,
\]
which is \eqref{eq:teles} with $C=e^{k_{0}^{q}}$.
\end{proof}

\begin{remark}\label{rem:telesconst}
In the setting of Lemma \ref{lem:telescope} let $p=p(\gamma)$ and $q=q(\gamma)$ satisfy $1<p<1/(\gamma-1)$ and $p\gamma<q<\gamma/(\gamma-1)$. By the proof of that lemma the constant in \eqref{eq:teles} can be taken to be $\exp(k^{q})$ for some $k\ge1$ with
\[
\log k\le K\Big(1+\log C_{0}+\log\frac{1}{\delta}+\log\Big(1+\frac{1}{a_{0}}\Big)\Big),
\]
where $K=K(\gamma)$. Moreover the choice of $k$ can be made so that $\exp(k^{q})$ is nondecreasing in $C_{0}$, in $\delta^{-1}$ and in $a_{0}^{-1}$.
\end{remark}

\begin{remark}\label{rem:telescope}
Lemma \ref{lem:telescope} is \cite[Proposition~2.2]{WZ26} extended to $0<\delta\le1$, the extension being needed because the exponent furnished by Corollary \ref{cor:onet} is smaller than one. The level set subdivision of the time axis replaces the geometric partition used in
\cite[proof of Theorem~1]{AEW}. A geometric partition has increments comparable to
$z^{-m}$, and the weights attached to it are summable only when the exponent of the
one-time inequality equals one, as in \cite[proof of Theorem~1]{AEW}. The increments of
the subdivision decay only polynomially, $\ell_k\ge c_0k^{-p}$, and the weight
$w_k=e^{-k^{q}}$ absorbs them precisely under the condition $p\gamma<q$.
\end{remark}

\subsection{Duality}\label{subsec:duality}
The last auxiliary result passes from the adjoint equation to the controlled one.
\begin{lemma}[Duality identity]\label{lem:duality}
Assume \eqref{eq:coef}--\eqref{eq:bound} and let all functions below be real-valued. Let $y$ be the solution of \eqref{eq:ctrl} with $y_{0}\in L^{2}_{\mathcal{F}_{0}}(\Omega;L^{2}(\T^{n}))$ and with $u\in L^{2}_{\mathbb{F}}(0,T;L^{2}(\T^{n}))$ vanishing outside $D$, and let $(z,Z)$ be the solution of \eqref{eq:bwd} with $z_{T}\in L^{2}_{\mathcal{F}_{T}}(\Omega;L^{2}(\T^{n}))$. Then
\begin{equation}\label{eq:duality}
\E\langle y(T),z_{T}\rangle-\E\langle y_{0},z(0)\rangle=\E\int_{0}^{T}\langle u(t),\chi_{D}z(t)\rangle\dd t .
\end{equation}
\end{lemma}

\begin{proof}
For $N\in\mathbb{N}$ let $H_{N}$ be the span in $L^{2}(\T^{n})$ of $\{e^{ik\cdot x}\ \text{with}\ \abs{k_{j}}\le N\ \text{for}\ j=1,\dots,n\}$ and let $P_{N}$ be the orthogonal projection of $L^{2}(\T^{n})$ onto $H_{N}$. Let $y_{N}$ be the $H_{N}$-valued solution of
\[
\dd y_{N}=\big(\Delta y_{N}-P_{N}(ay_{N})+P_{N}(\chi_{D}u)\big)\dd t+P_{N}(by_{N})\dd W(t),\qquad y_{N}(0)=P_{N}y_{0},
\]
and let $(z_{N},Z_{N})$ be the $H_{N}$-valued solution of
\[
\dd z_{N}=\big(-\Delta z_{N}+P_{N}(az_{N})-P_{N}(bZ_{N})\big)\dd t+Z_{N}\dd W(t),\qquad z_{N}(T)=P_{N}z_{T}.
\]
Both equations are finite dimensional and linear, with coefficients bounded in $L^{2}(0,T)$ in time, and hence are uniquely solvable. The energy estimate for the two difference equations, with the weights $\exp\big(\pm\int_{0}^{t}(2\alpha+2\beta^{2})\dd\tau\big)$ as in the proof of Lemma \ref{lem:scale}, gives
\[
y_{N}\to y\quad\text{and}\quad z_{N}\to z\quad\text{in }L^{2}_{\mathbb{F}}(\Omega;C([0,T];L^{2}(\T^{n}))),\qquad Z_{N}\to Z\quad\text{in }L^{2}_{\mathbb{F}}(0,T;L^{2}(\T^{n})).
\]
Since $P_{N}$ is self-adjoint and $y_{N},z_{N},Z_{N}$ take values in $H_{N}$, the projections may be omitted in the pairings below. It\^o's formula gives
\begin{equation*}
    \begin{split}
        \dd\langle y_{N},z_{N}\rangle=&\Big(\langle\Delta y_{N},z_{N}\rangle+\langle y_{N},-\Delta z_{N}\rangle+\langle\chi_{D}u,z_{N}\rangle-\langle ay_{N},z_{N}\rangle+\langle y_{N},az_{N}\rangle-\langle y_{N},bZ_{N}\rangle+\langle by_{N},Z_{N}\rangle\Big)\dd t\\
        &+\dd M_{N}
    \end{split}
\end{equation*}
with $M_{N}$ a continuous local martingale. Since $\Delta$ is self-adjoint and $a$ and $b$ are real-valued,
\[
\langle\Delta y_{N},z_{N}\rangle+\langle y_{N},-\Delta z_{N}\rangle=0,\qquad
\langle ay_{N},z_{N}\rangle=\langle y_{N},az_{N}\rangle,\qquad
\langle y_{N},bZ_{N}\rangle=\langle by_{N},Z_{N}\rangle,
\]
so that
\[
\dd\langle y_{N},z_{N}\rangle=\langle\chi_{D}u,z_{N}\rangle\dd t+\dd M_{N}.
\]
The local martingale has zero expectation by the stopping argument of Lemma \ref{lem:energy}. Hence
\[
\E\langle y_{N}(T),P_{N}z_{T}\rangle-\E\langle P_{N}y_{0},z_{N}(0)\rangle=\E\int_{0}^{T}\langle\chi_{D}u,z_{N}\rangle\dd t .
\]
Letting $N\to\infty$, and using the convergence above together with the strong convergence of $P_{N}$ to the identity, gives \eqref{eq:duality}, since $u=\chi_{D}u$.
\end{proof}

\section{Proof of the main results}\label{sec:proof}

\subsection{Proof of Theorem \ref{thm:obs}}\label{subsec:pfobs}

\begin{proof}
Let $z$ be a solution of \eqref{eq:bwd}. If $\E\norm{z(0)}^{2}=0$ there is nothing to prove, so assume $\E\norm{z(0)}^{2}>0$. Then Lemma \ref{lem:energy} gives $\E\norm{z(t)}^{2}>0$ and $h(t)>0$ for every $t\in[0,T]$. For $t\in(0,T)$ put
\[
g(t)=e^{\mu t/(2T)}\norm{z(t)}_{L^{2}_{\mathcal{F}_{t}}(\Omega;L^{2}(D_{t}))}.
\]

Step 1. We first select the times at which the spatial slice carries a fixed fraction of the measure. Put
\[
A_{D}=\Big\{t\in(0,T)\ \text{with}\ \abs{D_{t}}\ge\frac{d}{2T}\Big\},\qquad \nu=\frac{d}{2T\abs{\T^{n}}}.
\]
By Tonelli's theorem the functions $t\mapsto\abs{D_{t}}$ and $t\mapsto g(t)$ are measurable, hence $A_{D}$ is measurable and
\[
d\le\int_{0}^{T}\abs{D_{t}}\dd t\le\frac{d}{2}+\abs{\T^{n}}\abs{A_{D}} .
\]
Therefore $\abs{A_{D}}\ge d/(2\abs{\T^{n}})$ and $\abs{D_{t}}\ge\nu\abs{\T^{n}}$ for every $t\in A_{D}$. Moreover $g\in L^{1}(A_{D})$, since by the Cauchy--Schwarz inequality in time
\[
\int_{A_{D}}g(t)\dd t\le e^{\mu/2}T^{1/2}\Big(\int_{0}^{T}\E\norm{z(t)}_{L^{2}(\T^{n})}^{2}\dd t\Big)^{1/2}<+\infty .
\]

Step 2. Let $\vartheta\in(0,1)$ and $C_{\nu}\ge1$ be the constants of Corollary \ref{cor:onet} corresponding to this $\nu$. Let $t\in A_{D}$ and $s\in(t,T]$. By \eqref{eq:onet} with $E=D_{t}$,
\[
\norm{z(t)}_{L^{2}_{\mathcal{F}_{t}}(\Omega;L^{2}(\T^{n}))}\le C_{\nu}e^{C_{\nu}/(s-t)}
\norm{z(t)}^{\vartheta}_{L^{2}_{\mathcal{F}_{t}}(\Omega;L^{2}(D_{t}))}
\norm{z(s)}^{1-\vartheta}_{L^{2}_{\mathcal{F}_{s}}(\Omega;L^{2}(\T^{n}))}.
\]
By \eqref{eq:hcomp} one has
\[
h(t)\le e^{\mu/2}e^{\mu t/(2T)}\norm{z(t)}_{L^{2}_{\mathcal{F}_{t}}(\Omega;L^{2}(\T^{n}))},
\]
and, since $t<s\le T$ gives $t-s(1-\vartheta)<\vartheta s\le\vartheta T$,
\[
e^{\mu t/(2T)}\norm{z(s)}^{1-\vartheta}_{L^{2}_{\mathcal{F}_{s}}(\Omega;L^{2}(\T^{n}))}\le e^{\mu\vartheta/2}h(s)^{1-\vartheta}.
\]
By the definition of $g$ one has $\norm{z(t)}_{L^{2}_{\mathcal{F}_{t}}(\Omega;L^{2}(D_{t}))}=e^{-\mu t/(2T)}g(t)\le g(t)$. Since $e^{\mu/2+\mu\vartheta/2}\le e^{\mu}$, and since $\Lambda_{s,t}\ge1$ and $\Lambda_{s,t}\ge(s-t)^{-1}$ give $(s-t)^{-1}\le\Lambda_{s,t}^{3/2}$, one has $C_{\nu}/(s-t)\le C_{0}\Lambda_{s,t}^{3/2}$ with $C_{0}=C_{\nu}e^{\mu}$, and the last three displays give
\begin{equation}\label{eq:interp}
h(t)\le C_{0}\exp\big(C_{0}\Lambda_{s,t}^{3/2}\big)g(t)^{\vartheta}h(s)^{1-\vartheta}
\end{equation}
for every $t\in A_{D}$ and every $s\in(t,T]$.

Step 3. For $\tau\in[0,T]$ put $H(\tau)=h(T-\tau)$, $\tilde g(\tau)=g(T-\tau)$ and $B=T-A_{D}$. Then $H\in C([0,T])$ is positive and nonincreasing, $B\subset(0,T)$ is measurable with $\abs{B}=\abs{A_{D}}$, and $\tilde g\in L^{1}(B)$. Let $0\le\sigma<\tau\le T$ with $\tau\in B$ and put $t=T-\tau$, $s=T-\sigma$. Then $t\in A_{D}$, $s-t=\tau-\sigma$, and $\abs{D_{t}}\ge\nu\abs{\T^{n}}$ for every $\tau\in B$. Since $H$ is nonincreasing and $\Lambda_{s,t}=1+(\tau-\sigma)^{-1}$, \eqref{eq:interp} gives
\[
H(\tau)\le C_{0}\exp\Big(C_{0}\Big(1+(\tau-\sigma)^{-1}+\log\frac{H(\sigma)}{H(\tau)}\Big)^{3/2}\Big)\tilde g(\tau)^{\vartheta}H(\sigma)^{1-\vartheta}
\]
for all $0\le\sigma<\tau\le T$ and almost every $\tau\in B$. Since $d\le\abs{D}\le\abs{\T^{n}}T$, the number $a_{0}=d/(2\abs{\T^{n}})$ belongs to $(0,T)$, and $\abs{B}\ge a_{0}$ by Step 1. Hence Lemma \ref{lem:telescope} with $\gamma=3/2$, $\delta=\vartheta$, $A=B$, $h=H$ and $g=\tilde g$ yields
\begin{equation}\label{eq:obsH}
H(T)\le C_{\mathrm{tel}}\int_{B}\tilde g(\tau)\dd\tau=C_{\mathrm{tel}}\int_{A_{D}}g(t)\dd t,
\end{equation}
where $C_{\mathrm{tel}}=C(C_{0},3/2,\vartheta,T,a_{0})$.

Step 4. Since $H(T)=h(0)=\norm{z(0)}_{L^{2}_{\mathcal{F}_{0}}(\Omega;L^{2}(\T^{n}))}$, estimate \eqref{eq:obsH} gives
\[
\norm{z(0)}_{L^{2}_{\mathcal{F}_{0}}(\Omega;L^{2}(\T^{n}))}
\le C_{\mathrm{tel}}\int_{A_{D}}g(t)\dd t
\le C_{\mathrm{tel}}e^{\mu/2}\int_{0}^{T}\norm{z(t)}_{L^{2}_{\mathcal{F}_{t}}(\Omega;L^{2}(D_{t}))}\dd t ,
\]
which is \eqref{eq:obsL1} with $C_{obs}=C_{\mathrm{tel}}e^{\mu/2}$. Squaring and applying the Cauchy--Schwarz inequality in time together with Tonelli's theorem gives
\[
\E\norm{z(0)}^{2}
\le C_{\mathrm{tel}}^{2}e^{\mu}T\int_{0}^{T}\E\norm{z(t)}^{2}_{L^{2}(D_{t})}\dd t
= C_{obs}^{2}T\E\iint_{D}\abs{z(x,t)}^{2}\dd x\dd t ,
\]
so that \eqref{eq:obs} follows once \eqref{eq:obsL1} is available with a constant of the same form.

Step 5. It remains to convert the constants into the form stated. Recall $a_{0}=d/(2\abs{\T^{n}})$ and $\nu=d/(2T\abs{\T^{n}})$, hence $\nu^{-1}\le 2T\abs{\T^{n}}d^{-1}$ and $a_{0}^{-1}\le 2\abs{\T^{n}}d^{-1}$. Let $\kappa_{1},\kappa_{2},\kappa_{3}$ be the constants of Lemma \ref{lem:hilb} for $r=R/2$, so that $\vartheta^{-1}\le\kappa_{1}\big(1+\log(1/\nu)\big)$ and $C_{h}\le\kappa_{2}\nu^{-\kappa_{3}}$. By the definition of $C_{\nu}$ in the proof of Corollary \ref{cor:onet} there are $\kappa_{4}=\kappa_{4}(n,R)$ and $\kappa_{5}=\kappa_{5}(n,R,C_{1},T)$ such that
\[
\log C_{0}\le\kappa_{4}\log\frac{1}{\nu}+\kappa_{5}.
\]
By Remark \ref{rem:telesconst} with $\gamma=3/2$ and $\delta=\vartheta$ there is a constant $\kappa_{6}=\kappa_{6}(3/2)$ with
\[
\log\log C_{\mathrm{tel}}\le\kappa_{6}\Big(1+\log C_{0}+\log\frac{1}{\vartheta}+\log\Big(1+\frac{1}{a_{0}}\Big)\Big).
\]
Since $\log(1/\nu)=\log(1/d)+\log(2T\abs{\T^{n}})$, since $\log(1+\log(1/\nu))\le\log(1/\nu)$ gives $\log(1/\vartheta)\le\log\kappa_{1}+\log(1/\nu)$, and since $1\le\abs{\T^{n}}Td^{-1}$ gives $\log(1+1/a_{0})\le\log\big(\abs{\T^{n}}(T+2)\big)+\log(1/d)$, the last three displays give
\[
\log\log C_{\mathrm{tel}}\le\kappa_{7}\log\frac{1}{d}+\kappa_{8},
\qquad
\kappa_{7}=\kappa_{6}(\kappa_{4}+2),
\qquad
\kappa_{8}=\kappa_{8}(n,R,C_{1},T).
\]
Hence $C_{\mathrm{tel}}\le\exp\big(e^{\kappa_{8}}d^{-\kappa_{7}}\big)$. Since $C_{obs}=C_{\mathrm{tel}}e^{\mu/2}$ by Step 4 and since $d^{-\kappa_{7}}\le(1+d^{-1})^{\kappa_{7}}$ together with $1\le(1+d^{-1})^{\kappa_{7}}$, one has
\[
C_{obs}\le\exp\Big(\big(e^{\kappa_{8}}+\mu/2\big)(1+d^{-1})^{\kappa_{7}}\Big),
\qquad
C_{obs}^{2}T\le\exp\Big(\big(2e^{\kappa_{8}}+\mu+\abs{\log T}\big)(1+d^{-1})^{\kappa_{7}}\Big).
\]
Therefore \eqref{eq:obsL1} and \eqref{eq:obs} hold simultaneously with
\[
K=\kappa_{7},\qquad C=2e^{\kappa_{8}}+\mu+\abs{\log T},
\]
and this $C$ depends only on $n$, $R$, $C_{1}$ and $T$.
\end{proof}

\subsection{Proof of Theorem \ref{thm:ctrl}}\label{subsec:pfctrl}

\begin{proof}
In this proof all functions are real-valued, so that $\langle f,g\rangle=\int_{\T^{n}}fg\dd x$ for $f,g\in L^{2}(\T^{n})$. Let $C_{obs}$ be the constant of \eqref{eq:obsL1}.

Step 1. Put
\[
X=\big\{\chi_{D}z:(z,Z)\ \text{solves \eqref{eq:bwd} with}\ z_{T}\in L^{2}_{\mathcal{F}_{T}}(\Omega;L^{2}(\T^{n}))\big\}.
\]
Since \eqref{eq:bwd} is linear and its solution is unique for a given terminal datum, $X$ is a linear subspace of $L^{1}_{\mathbb{F}}(0,T;L^{2}(\Omega;L^{2}(\T^{n})))$. Define $L:X\to\R$ by $L(\chi_{D}z)=-\E\langle y_{0},z(0)\rangle$. If $\chi_{D}z=0$, then \eqref{eq:obsL1} gives $\norm{z(0)}_{L^{2}_{\mathcal{F}_{0}}(\Omega;L^{2}(\T^{n}))}=0$, so $L$ is well defined.

Step 2. By \eqref{eq:obsL1} and the Cauchy--Schwarz inequality,
\[
\abs{L(\chi_{D}z)}\le\norm{y_{0}}_{L^{2}_{\mathcal{F}_{0}}(\Omega;L^{2}(\T^{n}))}\norm{z(0)}_{L^{2}_{\mathcal{F}_{0}}(\Omega;L^{2}(\T^{n}))}
\le C_{obs}\norm{y_{0}}_{L^{2}_{\mathcal{F}_{0}}(\Omega;L^{2}(\T^{n}))}\norm{\chi_{D}z}_{L^{1}_{\mathbb{F}}(0,T;L^{2}(\Omega;L^{2}(\T^{n})))} ,
\]
so $L$ is bounded with norm at most $C_{obs}\norm{y_{0}}_{L^{2}_{\mathcal{F}_{0}}(\Omega;L^{2}(\T^{n}))}$. By the Hahn--Banach theorem $L$ extends to a bounded linear functional on $L^{1}_{\mathbb{F}}(0,T;L^{2}(\Omega;L^{2}(\T^{n})))$ without increase of norm, and by \cite[Theorem~2.73]{LZ} with $p=1$ and $q=2$ there is $u\in L^{\infty}_{\mathbb{F}}(0,T;L^{2}(\Omega;L^{2}(\T^{n})))$ such that
\[
\E\int_{0}^{T}\langle u(t),\chi_{D}z(t)\rangle\dd t=L(\chi_{D}z)
\]
for every $(z,Z)$ solving \eqref{eq:bwd} with $z_{T}\in L^{2}_{\mathcal{F}_{T}}(\Omega;L^{2}(\T^{n}))$, and
\[
\norm{u}_{L^{\infty}_{\mathbb{F}}(0,T;L^{2}(\Omega;L^{2}(\T^{n})))}\le C_{obs}\norm{y_{0}}_{L^{2}_{\mathcal{F}_{0}}(\Omega;L^{2}(\T^{n}))}.
\]
Since $\chi_{D}z=\chi_{D}(\chi_{D}z)$, replacing $u$ by $\chi_{D}u$ preserves the last identity and does not increase the norm, so we may assume that $u$ vanishes outside $D$. Since $L^{\infty}_{\mathbb{F}}(0,T;L^{2}(\Omega;L^{2}(\T^{n})))\subset L^{2}_{\mathbb{F}}(0,T;L^{2}(\T^{n}))$, the control $u$ is admissible for \eqref{eq:ctrl}.

Step 3. Let $y$ be the solution of \eqref{eq:ctrl} with this $u$ and with initial datum $y_{0}$. By Lemma \ref{lem:duality},
\[
\E\langle y(T),z_{T}\rangle-\E\langle y_{0},z(0)\rangle=\E\int_{0}^{T}\langle u(t),\chi_{D}z(t)\rangle\dd t
\]
for every $(z,Z)$ solving \eqref{eq:bwd} with $z_{T}\in L^{2}_{\mathcal{F}_{T}}(\Omega;L^{2}(\T^{n}))$. By the choice of $u$,
\[
\E\int_{0}^{T}\langle u(t),\chi_{D}z(t)\rangle\dd t=L(\chi_{D}z)=-\E\langle y_{0},z(0)\rangle .
\]
Hence $\E\langle y(T),z_{T}\rangle=0$ for every $z_{T}\in L^{2}_{\mathcal{F}_{T}}(\Omega;L^{2}(\T^{n}))$. Taking $z_{T}=y(T)$ gives $\E\norm{y(T)}^{2}=0$, so $y(T)=0$ almost surely, and \eqref{eq:ctrlcostLinf} is the bound of Step 2. Finally,
\[
\E\norm{u}^{2}_{L^{2}(0,T;L^{2}(\T^{n}))}=\int_{0}^{T}\E\norm{u(t)}^{2}\dd t\le T\norm{u}^{2}_{L^{\infty}_{\mathbb{F}}(0,T;L^{2}(\Omega;L^{2}(\T^{n})))},
\]
which together with \eqref{eq:ctrlcostLinf} and $C_{obs}^{2}T\le\exp\big(C(1+d^{-1})^{K}\big)$ gives \eqref{eq:ctrlcost}.
\end{proof}

\subsection{Proof of Proposition \ref{prop:orth}}\label{subsec:pforth}

\begin{proof}
By \eqref{eq:reach} and Lemma \ref{lem:duality} with $y_{0}=0$,
\begin{equation}\label{eq:orthduality}
\E\langle y_{u}(T),\eta\rangle=\E\int_{0}^{T}\langle u(t),\chi_{D}z_{\eta}(t)\rangle\dd t
\end{equation}
for every $u\in L^{2}_{\mathbb{F}}(0,T;L^{2}(\T^{n}))$ vanishing outside $D$.

Step 1. Assume that $\eta$ is orthogonal to $\mathcal{Y}$ and take $u=\chi_{D}z_{\eta}$ in \eqref{eq:orthduality}. Then $u\in L^{2}_{\mathbb{F}}(0,T;L^{2}(\T^{n}))$ and $u$ vanishes outside $D$, so the left hand side of \eqref{eq:orthduality} vanishes. By Tonelli's theorem the right hand side equals $\int_{0}^{T}\E\norm{z_{\eta}(t)}^{2}_{L^{2}(D_{t})}\dd t$, hence
\[
\E\norm{z_{\eta}(t)}^{2}_{L^{2}(D_{t})}=0\qquad\text{for almost every }t\in(0,T).
\]
Since $\abs{D_{t}}=0$ for $t\notin J$, this quantity vanishes for almost every $t\in J$.

Let $t\in J$ be such a time. By Lemma \ref{lem:scale} with $s=T$, the random variable $z_{\eta}(t)$ is the trace on $\T^{n}$ of a function $\mathcal Z_{t}$ holomorphic on $S^{n}_{3R/4}$ with values in $\mathcal H=L^{2}_{\mathcal{F}_{t}}(\Omega;\C)$, and
\[
M:=\sup_{\zeta\in S^{n}_{R/2}}\norm{\mathcal Z_{t}(\zeta)}_{\mathcal H}
\le C(n,R)\exp\Big(\frac{C_{1}^{2}}{2}+C_{1}T^{1/2}+\frac{9nR^{2}}{64(T-t)}\Big)\norm{\eta}_{L^{2}_{\mathcal{F}_{T}}(\Omega;L^{2}(\T^{n}))}<+\infty
\]
by \eqref{eq:analytic}. By Lemma \ref{lem:hilb} with $r=R/2$ and $E=D_{t}$ there are $C\ge1$ and $\vartheta\in(0,1)$, depending only on $n$, $R$ and $\abs{D_{t}}$, such that
\[
\norm{z_{\eta}(t)}_{L^{2}_{\mathcal{F}_{t}}(\Omega;L^{2}(\T^{n}))}
\le C\Big(\int_{D_{t}}\E\abs{z_{\eta}(t,x)}^{2}\dd x\Big)^{\vartheta/2}M^{1-\vartheta}=0 .
\]
Hence $z_{\eta}(t)=0$ in $L^{2}_{\mathcal{F}_{t}}(\Omega;L^{2}(\T^{n}))$ for almost every $t\in J$, which is \eqref{eq:orth}.

Step 2. Assume \eqref{eq:orth}. For $t\in(0,T)\setminus J$ one has $\abs{D_{t}}=0$ and hence $\E\norm{z_{\eta}(t)}^{2}_{L^{2}(D_{t})}=0$. Together with \eqref{eq:orth} this gives $\chi_{D}z_{\eta}=0$ in $L^{2}_{\mathbb{F}}(0,T;L^{2}(\T^{n}))$. By \eqref{eq:orthduality}, $\E\langle y_{u}(T),\eta\rangle=0$ for every $u\in L^{2}_{\mathbb{F}}(0,T;L^{2}(\T^{n}))$ vanishing outside $D$, that is, $\eta$ is orthogonal to $\mathcal{Y}$.
\end{proof}

\subsection{Proof of Theorem \ref{thm:approx}}\label{subsec:pfapprox}

\begin{proof}
By \eqref{eq:reach} and the linearity of \eqref{eq:ctrl}, the terminal values reachable from a given $y_{0}$ form the translate of $\mathcal{Y}$ by the terminal value of the solution of \eqref{eq:ctrl} with $u=0$. Hence \eqref{eq:ctrl} is approximately controllable at time $T$ if and only if $\mathcal{Y}$ is dense in $L^{2}_{\mathcal{F}_{T}}(\Omega;L^{2}(\T^{n}))$.

Sufficiency. Assume \eqref{eq:approxcond}. Let $\eta\in L^{2}_{\mathcal{F}_{T}}(\Omega;L^{2}(\T^{n}))$ be orthogonal to $\mathcal{Y}$ and let $(z_{\eta},Z_{\eta})$ be the solution of \eqref{eq:bwd} with $z_{T}=\eta$. By Proposition \ref{prop:orth}, one has $z_{\eta}(t)=0$ in $L^{2}_{\mathcal{F}_{t}}(\Omega;L^{2}(\T^{n}))$ for almost every $t\in J$, with $J$ given by \eqref{eq:obstime}.

For $\delta\in(0,T)$ one has
\[
\abs{D\cap(\T^{n}\times(T-\delta,T))}=\int_{T-\delta}^{T}\abs{D_{t}}\dd t ,
\]
so \eqref{eq:approxcond} gives $\abs{J\cap(T-\delta,T)}>0$ for every $\delta\in(0,T)$. Hence there is a sequence $t_{j}\to T$ such that $t_{j}\in J$ and $z_{\eta}(t_{j})=0$ for every $j$. Since $z_{\eta}\in L^{2}_{\mathbb{F}}(\Omega;C([0,T];L^{2}(\T^{n})))$, dominated convergence gives
\[
\lim_{j\to\infty}\E\norm{z_{\eta}(t_{j})-\eta}^{2}=0 ,
\]
hence $\eta=0$. Therefore $\mathcal{Y}^{\perp}=\{0\}$, and the linear subspace $\mathcal{Y}$ is dense.

Necessity. Let $\mathbf 1$ denote the constant function equal to $1$ on $\T^{n}$. Assume that \eqref{eq:approxcond} fails. Then there is $\delta_{0}\in(0,T)$ with
\[
\int_{T-\delta_{0}}^{T}\abs{D_{t}}\dd t=\abs{D\cap(\T^{n}\times(T-\delta_{0},T))}=0 ,
\]
hence $\abs{D_{t}}=0$ for almost every $t\in(T-\delta_{0},T)$. Fix $\delta\in(0,\delta_{0})$ and put
\[
J_{\delta}=\E\int_{T-\delta}^{T}\beta(r)^{2}\dd r .
\]
Since $\int_{0}^{T}\beta^{2}\le C_{1}^{2}$ almost surely by \eqref{eq:bound}, dominated convergence gives $J_{\delta}\to0$ as $\delta\to0$.

Step 1. Put $Z=\chi_{(T-\delta,T)}\mathbf 1$ and let $z$ be the solution on $[T-\delta,T]$ of
\begin{equation}\label{eq:sharpstate}
\dd z(r)=\big(-\Delta z(r)+a(r)z(r)-b(r)Z(r)\big)\dd r+Z(r)\dd W(r),\qquad z(T-\delta)=0 ,
\end{equation}
extended by $z=0$ on $[0,T-\delta]$. By \cite[Section 4]{LZ}, $z\in L^{2}_{\mathbb{F}}(\Omega;C([0,T];L^{2}(\T^{n})))\cap L^{2}_{\mathbb{F}}(0,T;H^{1}(\T^{n}))$ and $Z\in L^{2}_{\mathbb{F}}(0,T;L^{2}(\T^{n}))$.

Since $\langle-\Delta z,z\rangle=\norm{\nabla z}^{2}$, $\langle az,z\rangle\le\alpha\norm{z}^{2}$, $\abs{\langle bZ,z\rangle}\le\abs{\T^{n}}^{1/2}\beta\norm{z}$ and $\norm{Z}^{2}\le\abs{\T^{n}}$, It\^o's formula applied to $\norm{z}^{2}$ in \eqref{eq:sharpstate} gives, for $r\in[T-\delta,T]$,
\[
\norm{z(r)}^{2}\le\int_{T-\delta}^{r}\Big(\big(2\alpha+1\big)\norm{z}^{2}+\abs{\T^{n}}\big(1+\beta^{2}\big)\Big)\dd s+2\int_{T-\delta}^{r}\langle z,Z\rangle\dd W(s).
\]
Since $\int_{0}^{T}(2\alpha+1)\dd s\le2T^{1/2}C_{1}+T$ almost surely, the Burkholder--Davis--Gundy inequality and Gronwall's lemma give
\begin{equation}\label{eq:sharpenergy}
\E\sup_{T-\delta\le r\le T}\norm{z(r)}^{2}\le C\big(\delta+J_{\delta}\big),\qquad
\E\int_{T-\delta}^{T}\norm{\nabla z(r)}^{2}\dd r\le C\big(\delta+J_{\delta}\big),
\end{equation}
with $C=C(n,T,C_{1})$. The second estimate follows from the identity above before the gradient term is discarded, the stochastic integral there having zero mean because $\E\int_{T-\delta}^{T}\norm{z}^{2}\norm{Z}^{2}\dd r\le\abs{\T^{n}}\E\int_{T-\delta}^{T}\norm{z}^{2}\dd r<+\infty$. Hence $(z,Z)$ belongs to the solution space of \eqref{eq:bwd} and solves \eqref{eq:bwd} with terminal datum $\eta=z(T)$, so that $(z,Z)=(z_{\eta},Z_{\eta})$ by uniqueness. In particular $z_{\eta}(t)=0$ for every $t\in[0,T-\delta]$.

Step 2. Put $m(r)=\langle z(r),\mathbf 1\rangle=\int_{\T^{n}}z(r,x)\dd x$. Since $\langle\Delta z,\mathbf 1\rangle=\langle z,\Delta\mathbf 1\rangle=0$, \eqref{eq:sharpstate} gives
\[
m(T)=\int_{T-\delta}^{T}\Big(\int_{\T^{n}}az\dd x-\int_{\T^{n}}bZ\dd x\Big)\dd r+\abs{\T^{n}}\big(W(T)-W(T-\delta)\big).
\]
Since $\abs{\int_{\T^{n}}az\dd x}\le\abs{\T^{n}}^{1/2}\alpha\norm{z}$ and $\abs{\int_{\T^{n}}bZ\dd x}\le\abs{\T^{n}}\beta$, the Cauchy--Schwarz inequality in time, \eqref{eq:bound} and \eqref{eq:sharpenergy} yield
\[
\E\Big|\int_{T-\delta}^{T}\int_{\T^{n}}az\dd x\dd r\Big|^{2}\le\abs{\T^{n}}C_{1}^{2}\delta\,\E\sup_{T-\delta\le r\le T}\norm{z(r)}^{2}\le C\big(\delta^{2}+\delta J_{\delta}\big),
\]
\[
\E\Big|\int_{T-\delta}^{T}\int_{\T^{n}}bZ\dd x\dd r\Big|^{2}\le\abs{\T^{n}}^{2}\delta J_{\delta}.
\]
Since $\abs{x+y}^{2}\ge\tfrac12\abs{y}^{2}-\abs{x}^{2}$ and $\E\abs{W(T)-W(T-\delta)}^{2}=\delta$, the last three displays give
\[
\E\abs{m(T)}^{2}\ge\tfrac12\abs{\T^{n}}^{2}\delta-\delta\big(C\big(\delta+J_{\delta}\big)+\abs{\T^{n}}^{2}J_{\delta}\big).
\]
As $J_{\delta}\to0$ and $\delta\to0$, one may decrease $\delta$ so that $C(\delta+J_{\delta})+\abs{\T^{n}}^{2}J_{\delta}\le\tfrac14\abs{\T^{n}}^{2}$, and then $\E\abs{m(T)}^{2}\ge\tfrac14\abs{\T^{n}}^{2}\delta$. Since $\E\abs{m(T)}^{2}=\E\abs{\langle\eta,\mathbf 1\rangle}^{2}\le\abs{\T^{n}}\E\norm{\eta}^{2}$, it follows that $\eta\neq0$ in $L^{2}_{\mathcal{F}_{T}}(\Omega;L^{2}(\T^{n}))$.

Step 3. Since $\abs{D_{t}}=0$ for almost every $t\in(T-\delta,T)$ and $z_{\eta}(t)=0$ for every $t\in[0,T-\delta]$, one has $\chi_{D}z_{\eta}=0$ in $L^{2}_{\mathbb{F}}(0,T;L^{2}(\T^{n}))$. By Lemma \ref{lem:duality} with $y_{0}=0$,
\[
\E\langle y_{u}(T),\eta\rangle=\E\int_{0}^{T}\langle u(t),\chi_{D}z_{\eta}(t)\rangle\dd t=0
\]
for every $u\in L^{2}_{\mathbb{F}}(0,T;L^{2}(\T^{n}))$ vanishing outside $D$. Hence $\eta$ is orthogonal to $\mathcal{Y}$. Since $\eta\neq0$, the subspace $\mathcal{Y}$ is not dense.
\end{proof}

\section{Funding}
Not applicable.

\section{Declarations}
The authors have not disclosed any competing interests and data availability.

\bibliographystyle{abbrvnat}
\bibliography{ref}
\end{document}